\documentclass[11pt]{amsart}
\usepackage{amsmath,amsfonts,amsthm,amscd,amssymb,mathrsfs,amssymb}
\usepackage{mathrsfs}
\usepackage{graphicx}
\usepackage{wrapfig}
\usepackage[all]{xy}
\newtheorem{theorem}{Theorem}

\newtheorem{proposition}[theorem]{Proposition}
\newtheorem{lemma}[theorem]{Lemma}

\newtheorem{remark}[theorem]{Remark}
\newcommand{\aaa}{\alpha}

\newcommand{\lmd}{\lambda}
\newcommand{\Lmd}{\Lambda}
\newcommand{\eee}{\epsilon}
\newcommand{\CP}{\mathbb{CP}}
\newcommand{\CC}{\mathbb{C}}

\newcommand{\ZZ}{\mathbb{Z}}

\newcommand{\ol}{\overline}
\newcommand{\lra}{\longrightarrow}

\newcommand{\proofend}{\hfill$\square$}

\newcommand{\Ker}{{\rm{Ker}}}
\newcommand{\Bs}{{\rm{Bs}}}

\newcommand{\vsp}{\vspace{3mm}}

\numberwithin{equation}{section}
\numberwithin{theorem}{section}

\begin{document}
\bibliographystyle{alpha} 
\title{Double solid twistor spaces II: general case}
\author{Nobuhiro Honda}
\address{Mathematical Institute, Tohoku University,
Sendai, Miyagi, Japan}
\email{honda@math.tohoku.ac.jp}
\begin{abstract}
In this paper we investigate Moishezon twistor spaces
which have a structure of double covering over 
a very simple rational threefold. 
These spaces can be regarded as a direct generalization of 
the twistor spaces studied in the papers \cite{P92,KK92}
to the case of arbitrary signature.
In particular, the branch divisor of the double covering is a cut of the rational threefold by
a single quartic hypersurface.
A defining equation of the hypersurface is determined in an explicit form.
We also show that these twistor spaces interpolate
LeBrun twistor spaces and  the twistor spaces constructed 
in \cite{HonLBl}.

\end{abstract}
\maketitle
\setcounter{tocdepth}{1}
\vspace{-5mm}


\section{Introduction}
In an influential paper \cite{Hi81},
Hitchin initiated a systematic study of compact twistor spaces
by algebro-geometric means.
In particular, by investigating the half-anticanonical system of twistor spaces, he showed that if a compact twistor space admits a K\"ahler metric, it must be one of the two standard twistor spaces, the projective space
or the flag manifold.
This direction of research was succeeded by 
Poon \cite{P86,P92} and Kreussler-Kurke \cite{KK92} to determine structure of twistor spaces over the
connected sum of two or three complex projective planes, by means of the same system on the twistor spaces.
Also, LeBrun \cite{LB91} and Campana-Kreussler \cite{CK98} utilized the same system
to investigate particular Moishezon twistor spaces over the
connected sum of any number of complex projective planes.

However, for other twistor spaces, it was evident that 
the half-anticanonical system 
no longer brings enough information for analyzing  structure of the spaces, because the system is at most a pencil.
Therefore  multiples of the half-anticanonical system have been used in order to find out
and explore new Moishezon twistor spaces.
In these analyses the first essential part is always
to find a pluri-half-anticanonical system which is 
not composed with the 
half-anticanonical system. 
Once this is established, by investigating  
the rational map associated to the multiple  system,
we can make a traditional analysis in algebraic geometry
to gain detailed structure of 
the twistor spaces.

In the paper \cite{HonDSn1}, we pursued such a direction and found a series of Moishezon twistor spaces on $n\CP^2$, the connected sum of $n$ complex projective planes where $n$ being arbitrary
with $n\ge 4$, such that the $(n-2)$-th power of the half-anticanonical system induces a rational map which is two-to-one over the image.
This image is a scroll of 2-planes over a rational normal curve, which is canonically embedded in 
$\CP^n$, and
the branch divisor is a cut of the scroll by a single quartic hypersurface.
Further, a defining equation of the quartic hypersurface was 
determined. 

While these twistor spaces on $n\CP^2$
 can be regarded as a generalization of the twistor spaces
on $3\CP^2$
studied by Poon \cite{P92} and Kreussler-Kurke \cite{KK92}, from  detail investigation in the case of $4\CP^2$ \cite{HonDS4_1,
HonDS4_2}, it is strongly expected  that 
the twistor spaces in \cite{HonDSn1} are specialization
of more general twistor spaces which also have a double covering structure over the same scroll. 
A purpose of this paper is to show that this is really the case, and determine the defining equation
of the quartic hypersurface which 
cuts out the branch divisor of the double covering.
However, for the actual analysis  totally new method is required,
because  unlike the ones
in \cite{HonDSn1}, the present twistor spaces
do not have an effective $\CC^*$-action.


In Section \ref{ss:S} we  construct a rational surface
$S$ 
and  investigate its pluri-anticanonical systems.
We are concerned with a twistor space on $n\CP^2$
which contains this surface $S$ as a member of the
half-anticanonical system $|F|$.
%
We devote most of Section \ref{ss:ps} to prove that  the 
multiple system $|(n-2)F|$ of the twistor space is not composed with the system $|F|$ (which is just a pencil).
The basic idea for proving this is simple to the effect that
we pick up distinct $(n-2)$ general members of the pencil
$|F|$ and  look at the restriction of the system $|(n-2)F|$ to the union
of these members.
However, this does not work well in this 
primitive form and we need to 
blowup the twistor space at the base curve of $|F|$.
This makes the $(n-2)$ divisors disjoint, and further by letting the
exceptional divisor of the blowup to be included in
the restriction, we obtain a crucial vanishing  of 
cohomology groups (Proposition \ref{prop:mainvan}).
This reduces the computations for the multiple system
to those on a divisor of smooth normal crossing.
The computations on the last divisor work 
very effectively, and 
we can finally show that the system $|(n-2)F|$ is 
$n$-dimensional as a linear system (Proposition \ref{prop:Z}).
Once this is proved, it is not difficult to 
show that 
 the image of the rational map associated
 to the multiple system is a scroll of planes
over a rational normal curve in $\CP^{n-2}$,
and that the map is two-to-one over the scroll,
whose branch divisor is a cut of the scroll by a quartic hypersurface
(Proposition \ref{prop:Z2}).

Section \ref{s:red} is employed to show that 
there exist two special  reducible members
of the system $|(n-2)F|$ which consist of two irreducible components.
These two divisors play a significant role for
obtaining the defining equation of the quartic
hypersurface.
The Chern classes of the irreducible components 
of these reducible members are
presented
in  explicit forms (Proposition \ref{prop:red}).
The method for proving this existence is similar to the method
in Section \ref{s:ps}, but as a restriction we have to take degree-one divisors
instead of the divisors in $|F|$,
which makes the computations  considerably heavier and 
much more subtle.
But again by making the exceptional divisor
of the same blowup
to be included in the restriction, we are able to 
obtain a critical vanishing result (Proposition \ref{prop:mainvan3}).
Then after long computations over the degree-one divisors
and the exceptional divisors,
we finally obtain the desired existence
of the special reducible members of 
the system $|(n-2)F|$.

Like most other non-projective Moishezon twistor spaces,
the multiple system $|(n-2)F|$ has non-empty base locus.
In Section \ref{s:deg1} we provide
a complete elimination of the base locus,
by a succession of explicit blowups
(Proposition \ref{prop:elim}).
While the base locus of the system consists of just strings of 
smooth rational curves, 
for a complete elimination a number of blowing-up is required.
As a consequence, we obtain information
about the image of particular twistor lines
and degree-one divisors under the rational
map associated to $|(n-2)F|$.

In the final section, assembling all the results
 in Sections \ref{s:ps}, \ref{s:red} and \ref{s:deg1},
we determine a defining equation of the quartic hypersurface
which cuts out the branch divisor on the scroll.
In Section \ref{ss:dc}, by utilizing the reducible members obtained
in Section \ref{s:red}, we find special curves on the branch divisor and  
 prove the existence of a quadratic 
hypersurface in $\CP^n$ which contains all these
special  curves.
In Section \ref{s:de} we prove the main result which determines
the defining equation of the quartic hypersurface
(Theorem \ref{thm:main}).
The argument in the proof is mostly algebraic, and is an improvement of the proof 
given in \cite[Theorem 4.5]{HonDS4_1}.
Finally in Section \ref{ss:moduli} we first compute
the dimension of the moduli space of the present twistor spaces.
Next we discuss some global structure of the moduli space.
In particular, we see that the present moduli space
can be partially compactified (completed) by attaching some part (a stratum)
of the moduli spaces of LeBrun twistor spaces.
We also mention that 
a stratum of of the moduli space of  the twistor 
spaces constructed in \cite{HonLBl} is also naturally attached 
to give a partial compactification of the moduli space
of the present twistor spaces.

\vsp
\noindent
{\bf Notation.}
The letter $F$ always denotes the canonical square root
of the anticanonical line bundle over a twistor space.
The degree of a divisor in a twistor space means
the intersection number with a twistor line.
Linear equivalence between divisors are often denoted by `$\sim$'.
If $Y$ is a subvariety of a complex manifold $X$,
and if $\mathscr S$ is a sheaf over $X$,
we often write $H^q(Y,\mathscr S)$ to mean $H^q(Y,\mathscr S|_Y)$.
We also write $h^q(X,\mathscr S)$ 
for $\dim H^q(X,\mathscr S)$.
An $(a,b)$-curve on the product surface $\CP^1\times
\CP^1$ means a curve of bidegree $(a,b)$.
For a vector space $V$, $S^kV$ denotes the $k$-th symmetric
product.
We often identify a divisor with the associated line bundle.
Usually we use the same letter for an 
analytic subspace in a complex space
and its strict transform into a blowup.

\section{Analysis of the pluri-half-anticanonical system}
\label{s:ps}

\subsection{Construction of a rational surface}
\label{ss:S}
First we  construct a rational surface $S$ which will be 
contained in the twistor spaces as a real irreducible member of the system $|F|$.

Let $n\ge 4$ be any fixed integer.
First consider the product surface $\CP^1\times\CP^1$ and define
a real structure on it as the product of the complex conjugation and the anti-podal map.
Fix a real reducible $(2,2)$-curve consisting of four irreducible components.
We write it as $C_1+C_2+\ol C_1 +\ol C_2$, where $C_1$ and $\ol C_1$ are
$(1,0)$-curves and $C_2$ and $\ol C_2$ are $(0,1)$-curves.
Next we choose any two points on $C_1$ which are not on  $C_2\cup\ol C_2$ (i.e.\,not on the corners).
We also choose any one point on $C_2$ which is not on $C_1\cup \ol C_1$.
By taking the conjugation by the real structure of these 3 points, we obtain $6$ points on the $(2,2)$-curve.
Let $S_0\to \CP^1\times \CP^1$ be the blowup at these 6 points.
$S_0$ is equipped with a  real structure lifted from $\CP^1\times\CP^1$.

Next we blowup $S_0$ at the two  points 
$C_2\cap \ol C_1$ and $\ol C_2\cap C_1$ $(n-3)$ times respectively,
where if $n>4$ the blowup is always done in the direction of $\ol C_1$
and $C_1$ respectively.
(This means that blown-up point is always the intersection point of $C_1$ with the exceptional curves of the last blowup,
and similar for the conjugate point.)
Let $S\to S_0$ be the resulting birational morphism.
Since we have blown-up $2n$ times in total,
we obtain $K_S^2=8-2n$.
Let 
\begin{align}\label{cycle1}
C:= C_1+C_2+\cdots+C_{n-1}+\ol C_1+\ol C_2+\cdots
+\ol C_{n-1}
\end{align}
be the unique anticanonical curve on $S$ arranged
in a cyclic order.
From the construction it is immediate to see that the self-intersection numbers of the components $C_1,C_2,\cdots,C_{n-1}$ in $S$ are 
respectively given by
\begin{align}\label{B1}
1-n, {\overbrace{-2,-2,\cdots,-2}^{n-3}},-1.
\end{align}
These intersection numbers are of fundamental importance 
throughout this paper.
The original real structure 
on $\CP^1\times\CP^1$  naturally lifts to the surface $S$, under which
 the cycle $C$ is real.
We note that in this construction of the surface $S$,
all  freedom is in choosing two points on $C_1$
and one point on $C_2$ in the construction of the surface $S_0$,
and there is no freedom in all the remaining blowups $S\to S_0$.

\begin{proposition}\label{prop:S}
The pluri-anticanonical systems of the surface $S$ enjoy the following properties:
\begin{enumerate}
\item[(i)]
If $0\le m< n-2$, then $h^0(mK_S^{-1}) = 1$.
\item[(ii)] 
$h^0((n-2)K_S^{-1}) = 3$.
\item[(iii)] The fixed components of the system \,$|(n-2)K_S^{-1}|$ is the curve
\begin{align}
\label{bc0}
(n-3)(C_1+\ol C_1) + \sum_{i=2}^{n-2}(n-1-i)(C_i+ \ol C_i).
\end{align}
\item[(iv)]
After removing this curve, the system \,$|(n-2)K_S^{-1}|$  is base point free.
\item[(v)]
If $\phi:S\to \CP^2$ denotes the morphism associated to the last system,
$\phi$ is of degree two and 
the branch divisor is a quartic curve.

\end{enumerate}
\end{proposition}

\proof
(i) and (iii) can be proved by computing intersection numbers, and we omit the detail.
For (ii), if 
we define a line bundle $L$ over $S$ 
as $(n-2)K_S^{-1}$ minus the curve \eqref{bc0},
then
for any $i\neq 2$, the intersection number 
$(L, C_i)_S$ can be seen to be zero.
Hence we have an exact sequence
\begin{align}\label{ses1}
0 \lra L - (C-C_2-\ol C_2)
\lra L \lra \mathscr O_{C-C_2-\ol C_2}\lra 0.
\end{align}
If we write $L'$ for the first non-trivial term of this sequence,
by Riemann-Roch formula we readily obtain $\chi ( L' ) = 1$.
From intersection numbers we also get $h^0(L')=1$.
As the line bundle $L'$ has a non-zero section,
we also have $H^2 (L' ) = 0$.
Hence we obtain $H^1(L')=0$.
Therefore noting that the curve $C-C_2-\ol C_2$ consists of two connected components,
from \eqref{ses1}, we get $h^0(L) = 3$, and we obtain (ii).
For (iv) we readily have $(L,L)_S=2$ and $(L,C_2)_S=1$.
Also, from the cohomology exact sequence of \eqref{ses1} 
we obtain that $\Bs\,|L|$ is disjoint from the curve $C-C_2-\ol C_2$.
These mean that if $\Bs\,|L|\neq\emptyset$ then $\Bs\,|L|$ consists of 
two points, one of which is on $C_2$  and the other is on $\ol C_2$.
Moreover these two points are not
on the complement $C-C_2-\ol C_2$.
Let $\nu:S'\to S$ be the blowup at these two points, $E$ and $\ol E$
the exceptional curves,
and put $L''= \nu^*L- E-\ol E$.
Then we have $(L'',L'')_{S'}=(L,L)_S-2=0$.
As $(L, C_2)_S=1$ we also obtain that $\Bs\,|L''|=\emptyset$.
Hence the rational map $\psi:S'\to \CP^2$ 
associated to $|L''|$ is a morphism whose
image is a curve.
On the other hand, since $(\nu^*L- E-\ol E, E)_{S'} = 1$,
the image $\psi(E)$ must be a line.
Hence $\psi(S')$ is a line.
So $\phi(S)$ is also a line, which is a contradiction.
Therefore we obtain $\Bs\,|L|=\emptyset$, meaning (iv).
For (v), since $(L,L)_S=2$, the morphism $\phi:S\to\CP^2$ is of degree two. 
Also
 the arithmetic genus of $L$ is easily seen to be
one.
Hence the branch curve must be a quartic curve.
\endproof

\vspace{3mm}
We have more detail about the double covering map $\phi:S\to\CP^2$.
This would be useful for 
understanding singularities of the branch divisor of the double covering map 
from the twistor space that will be obtained at the end of this section:
\begin{proposition}
\label{prop:sing2}
Let $\phi:S\to \CP^2$ be the degree two morphism as in Proposition \ref{prop:S} (v).
Then $\phi$ maps the two curves $C_2$ and $\ol C_2$ to an identical line isomorphically,
and maps two connected curves $C_3\cup C_4\cup\cdots\cup C_{n-1}\cup \ol C_1$ and $\ol C_3\cup \ol C_4\cup\cdots\cup \ol C_{n-1}\cup  C_1$ to points
on the line.
Moreover, the branch quartic curve of $\phi$ has two ordinary nodes at these two points.
\end{proposition}

We do not give a proof for this proposition, since it can be shown in a standard way.

\subsection{Pluri-half-anticanonical systems of the twistor spaces}
\label{ss:ps}
We still fix any integer $n\ge 4$ and let $S$ be the surface 
obtained from $\CP^1\times\CP^1$ by blowing up $2n$ times
as in Section \ref{ss:S}.
Next let $Z$ be a twistor space on $n\CP^2$ and suppose
that $Z$ contains the surface $S$ as a real
member of the system $|F|$.
As before let $C$ be the unique anticanonical curve
\eqref{cycle1} on $S$.
Since $F|_S\simeq K_S^{-1}$, in view of Proposition \ref{prop:S},
it is tempting to expect that 
the restriction map
\begin{align}\label{rest1}
H^0(Z, (n-2)F) \lra H^0 (S, (n-2)F)\simeq H^0((n-2)K_S^{-1})
\end{align}
 is  surjective.
However, {\em for any $n> 4$, there is an example of Moishezon twistor space $Z$ on $n\CP^2$, such that $Z$ has a real smooth $S\in |F|$ whose $|(n-2)K_S^{-1}|$ fulfills the properties (i)--(v) (with some minor modification for the explicit form of the fixed component \eqref{bc0}) but nevertheless 
the restriction map \eqref{rest1} is {\em not} surjective.}
Thus  validity of the above expectation is very subtle.

In spite of this, for the present twistor spaces,
we have the following. 
\begin{proposition}\label{prop:Z}
For the pluri-half-anticanonical systems of the twistor space $Z$ on $n\CP^2$, we have the following.
\begin{enumerate}
\item[(i)] $h^0(F)=2$, and \,$\Bs\,|F| = C$.
\item[(ii)] If $m<n-2$, then $H^0(mF)=S^m H^0(F)$, so that
\,$h^0(mF) = m+1$.
\item[(iii)] 
$h^0((n-2)F) = n+1$.
In particular, the system $|(n-2)F|$ is not composed 
with the pencil \,$|F|$.
\end{enumerate}
\end{proposition}

The assertions (i) and (ii) can be shown in a standard way
by using Proposition \ref{prop:S}, and we omit a proof.
For the rest of this section  we give a proof of (iii).

For this, we first make it clear about reducible members
of the pencil $|F|$.
This issue is also now standard and we omit a proof.
For each $1\le i \le n-1$ let $L_i$ be the twistor line
through the point $C_i\cap C_{i+1}$
(i.e.\,a corner of the cycle $C$), where we read
$C_{n}=\ol C_1$. 
Then there exists a unique reducible member of $|F|$ which contains
the twistor line $L_i$.
This member consists of two irreducible components, and
their intersection  is exactly $L_i$.
We write this member of $|F|$ as $S_i^+ + S_i^-$, where 
we make distinction of the two components by promising that $S_i^-$ contains the curve $C_1$.

\begin{figure}
\includegraphics{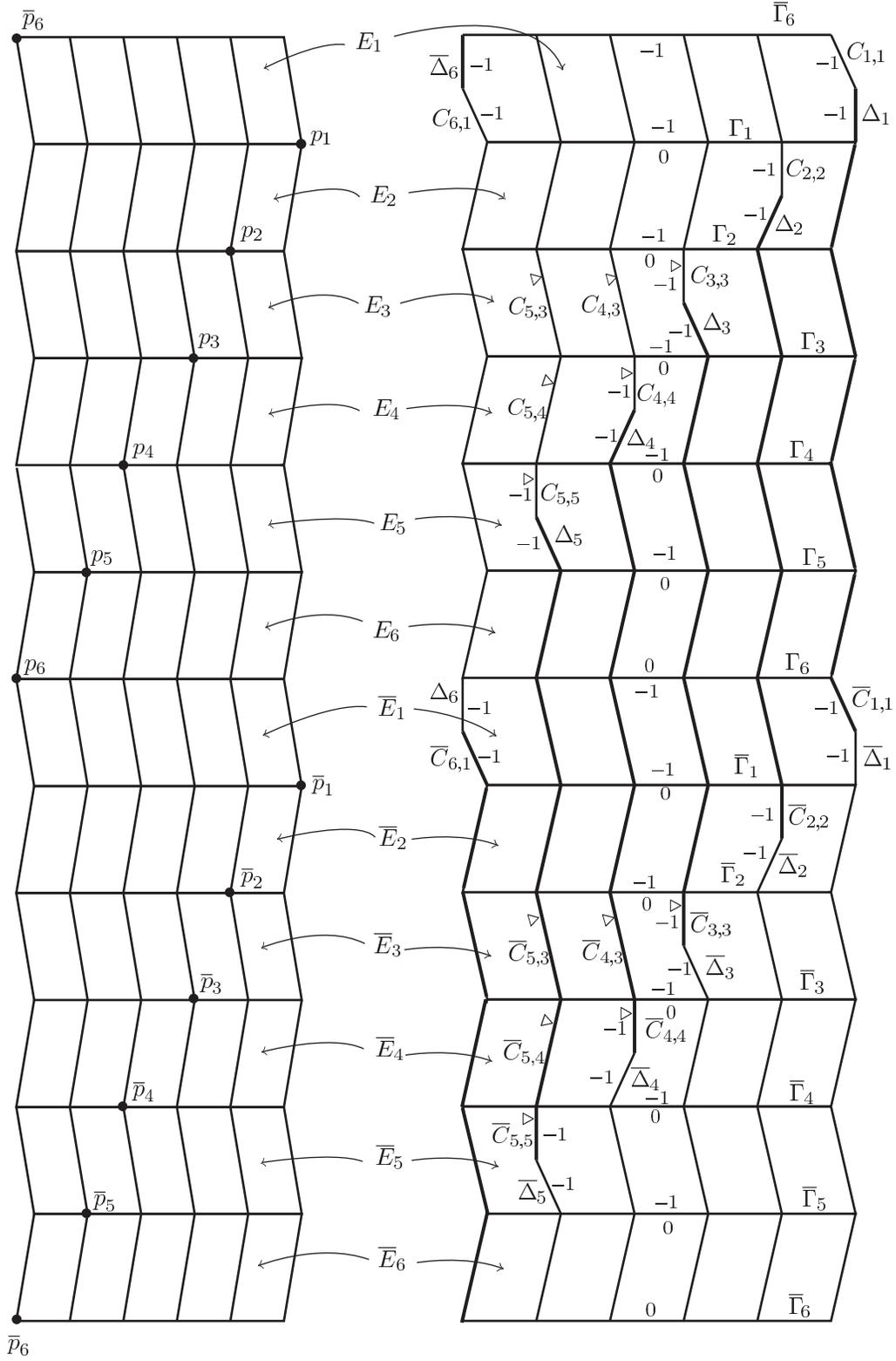}
\caption{The exceptional divisor of $\hat Z\to Z$
(left) and that of $\mu:Z_1\to Z$ (right)
in the case $n=7$.
In the right picture the intersection with $S_i^+$ 
($1\le i\le 6$)
is indicated by bold segments.
The segment with a small triangle represents
a base curve of $|\mathscr L_1|$ (Section \ref{s:deg1}). }
\label{fig:L1}
\end{figure}

Let $f:Z\to \CP^1$ be the rational map associated to the 
pencil $|F|$. 
The indeterminacy locus of $f$ is exactly the cycle $C$
by Proposition \ref{prop:Z} (i).
Let $\hat Z\to Z$ be the blowup of $Z$ at $C$, and $E_i$ and $\ol E_i$
$(1\le i\le n-1)$ the exceptional divisors over the components $C_i$ and $\ol C_i$ respectively.
The composition $\hat Z\to Z\to \CP^1$ is a morphism,
and has exactly $(n-1)$ reducible fibers, for which 
we still write $S_i^+\cup S_i^-$.
From the fact that $C$ constitutes a cycle,
it follows that every components $E_i$ and $\ol E_i$ are isomorphic
to $\CP^1\times\CP^1$, and on each of these components the morphism $\hat Z\to \CP^1$ 
coincides with a projection to one of the two factors.
For simplicity of notation we write the total sum of
the exceptional divisors as
$$
E:=\sum_{i=1}^{n-1}E_i + \sum_{i=1}^{n-1} \ol E_i.
$$
As the curve $C$ forms the cycle,
the divisor $E$ constitutes a `cylinder'
(see the left picture in Figure \ref{fig:L1}).
If we  denote the strict transform of the twistor line
$L_i$ ($1\le i\le n-1$) into $\hat Z$ by the same letter,
the intersection $L_i\cap E$ consists of two points,
and these are ordinary double points
(ODP-s for short) of 
the variety $\hat Z$.
(So $\hat Z$ has $2(n-1)$ ODP-s in total.)
For each $i$,
one of these two ODP-s is shared by
the four divisors 
$S_i^+, S_i^-, E_i$ and $E_{i+1}$,
and the other point is shared by
the divisors
$S_i^+,  S_i^-,  \ol E_i$ and $ \ol E_{i+1}$,
where if $i=n-1$ we read 
$E_n= \ol E_1$ and $\ol E_n= E_1$.
We denote by $p_i$ and $\ol p_i$ for the former and the latter
ODP respectively.
(These are indicated in the left picture in Figure \ref{fig:L1}.)

For each of these ODP-s there are two ways of small resolutions.
In order for later computations
to be transparent, we choose  small resolutions for  these in the following way.
When the index $i$ satisfies $1\le i\le n-2$, 
at the point $p_i$,
we take the small resolution which blows up
the pair  $\{S_i^+,E_i\}$.
When $i=n-1$, 
at the point $p_{n-1}$,
we take the small resolution which blows up
the alternative pair $\{S^-_{n-1},\ol E_1\}$.
For the conjugate point $\ol p_i$, $1\le i\le n-1$, we take the
small resolution which is determined from that of $p_i$
by the real structure.

Let $Z_1\to \hat Z$ be the birational morphism obtained by taking all these  small resolutions
simultaneously,
and $\Delta_i$ and $\ol{\Delta}_i$ ($1\le i\le n-1$) the 
exceptional curves of the points $p_i$ and $\ol p_i$ 
respectively.
We define a morphism $\mu:Z_1\to Z$ to be the composition $Z_1\to \hat Z\to Z$.
Also we write $f_1:Z_1\to \CP^1$ for the composition morphism
$Z_1\to \hat Z\to \CP^1$. 
This is precisely the rational map
associated with the pencil $|\mu^* F|$.
Of course $Z_1$ is non-singular
and equipped with a natural real structure.
Under the small resolution $Z_1\to \hat Z$,
each of the exceptional divisors $E_i$ 
(and $\ol E_i$) receives the following effect:
\begin{itemize}
\item the divisor $E_1\subset \hat Z$ is blownup at the two points $p_1$ and $\ol p_{n-1}$,
and the two curves $\Delta_1$ and $\ol{\Delta}_{n-1}$
are inserted as the exceptional curves,
\item
when $1<i<n-1$, the divisor $E_i\subset \hat Z$ is blownup at 
one point $p_i$, and the curve $\Delta_i$ is inserted as the exceptional curve,
\item
the divisor $E_{n-1}\subset \hat Z$ remains unchanged,
and hence the strict transform $E_{n-1} \subset Z_1$
is still biholomorphic to $\CP^1\times\CP^1$.
\end{itemize}
We illustrate these changes in Figure \ref{fig:L1}
in the case $n=7$.
The manifold $Z_1$ is the main stage for our computations,
in which these exceptional divisors play a principal role.

We use the same letters to mean the strict transforms 
into $Z_1$ of the
above divisors in $\hat Z$ and $Z$.
From Proposition \ref{prop:S} (iii), the pulled-back system 
$|\mu^*((n-2)F)|$ has the following divisor 
as fixed components at least:
\begin{align}\label{fc2}
(n-3)(E_1+\ol E_1) + \sum_{i=2}^{ n-2}(n-1-i)(E_i+ \ol E_i).
\end{align}
So we define a line bundle $\mathscr L_1$ over $Z_1$ by
\begin{align}\label{L}
\mathscr L_1:= \mu^*((n-2)F) 
- (n-3)(E_1+\ol E_1) - \sum_{i=2}^{n-2}(n-1-i)(E_i+ \ol E_i).
\end{align}
This is the line bundle we actually need to investigate.
Note that this line bundle is also real.
Here we outline our our method for obtaining $h^0(\mathscr L_1)$.

\begin{itemize}
\item[(i)]
We choose any distinct non-singular members $S_1,S_2,\cdots,S_{n-2}$ 
of the original pencil $|F|$. 
We use the same letters to mean the strict transforms of these divisors into $Z_1$.
\item[(ii)]
We restrict the line bundle $\mathscr L_1$ to the divisor
\begin{align}\label{restdiv}
(S_1\sqcup S_2\sqcup \cdots\sqcup S_{n-2}) \,\cup \,E,
\end{align}
where $E$ is the total sum of the exceptional divisors 
of the birational morphism $\mu:Z_1\to Z$.
The point here is that
unlike the union
$S_1\cup S_2\cup \cdots\cup S_{n-2}$ in
the original space $Z$,
 the divisor \eqref{restdiv} in $Z_1$ is clearly 
{\em smooth normal crossing}.
\item[(iii)] We prove that the restriction map $H^0(Z_1,\mathscr L_1)
\to H^0 ((S_1\sqcup \cdots\sqcup S_{n-2}) \,\cup \,E,\mathscr L_1)$ is surjective, 
by showing that the cohomology group $H^1$ of the kernel line bundle vanishes.
\item[(iv)]
 Here comes another point: all the restrictions of 
$\mathscr L_1$ to the divisors 
 $S_i$ and $E$ are explicitly computable.
Then thanks to the fact that 
the divisor \eqref{restdiv} is smooth normal crossing,
 it is possible to compute the space
 $H^0 ((S_1\sqcup \cdots\sqcup S_{n-2}) \,\cup \,E,\mathscr L_1)$.
\end{itemize}

\begin{remark}
{\em
Since the above argument might look a  bit complicated,
it should be explained why we consider the restriction
of the line bundle $\mathscr L_1$ to the divisor \eqref{restdiv}
for computing $h^0((n-2)F)$. 
By choosing the divisors $S_1,\cdots, S_{n-2}\in |F|$ as
above, we have the exact sequence
$
0 \to \mathscr O_Z \to (n-2)F \to (n-2)F\,|_{S_1\cup\cdots\cup S_{n-2}}\to 0
$
on $Z$.
Further we have $H^1(\mathscr O_Z)=0$
as $Z$ is simply connected.
So we can determine $h^0((n-2)F)$ if we could compute
$h^0((n-2)F\,|_{S_1\cup\cdots\cup S_{n-2}})$.
As $(n-2)F|_{S_i}\simeq (n-2)K_{S_i}^{-1}$
for each $i$ and as we know $h^0((n-2)K_{S_i}^{-1})=3$ by
Proposition \ref{prop:S}, one might think  it possible to compute
$h^0((n-2)F\,|_{S_1\cup\cdots\cup S_{n-2}})$.
However, this seems to be impossible due to the fact that the union
$S_1\cup \cdots\cup S_{n-2}$ is {\em not} smooth normal crossing,
when $n>4$.
This situation can be resolved by blowing up  the base curve $C$. 
But it is still impossible to compute $h^0(\mathscr L_1)$ if we just restrict the line bundle 
$\mathscr L_1$ to the  disjoint union $S_1\sqcup\cdots \sqcup S_{n-2}$ of the strict transforms, because this time
 we cannot expect the restriction map 
to be surjective (at the level of sections). 
Hence we make the divisor $E$ to be
included in the restriction.
As we see below, this method works very effectively.
}
\end{remark}
Beginning the actual computations for  $h^0(\mathscr L_1)$, 
we first compute the restriction
of $\mathscr L_1$ to the exceptional divisor $E_i$.
For this, we first define curves on $Z_1$ by
\begin{align}\label{cij}
C_{i,j} := (S_i^+\cup S_i^-) \cap E_j,
\quad
i,j\in \{1,2,\cdots,n-1\}.
\end{align}
This curve is naturally identified with 
the rational curve $C_j$ in $Z$ through the birational morphism $\mu$, and the first index $i$ indicates in which degree-one
divisor  the curve is contained.
We also define other curves on $Z_1$ by
\begin{align}\label{gi}
\Gamma_i := E_i\cap E_{i+1},
\quad
1\le i\le n-1,
\end{align}
where we read $E_n=\ol E_1$ when $i=n-1$.
See the right picture in Figure \ref{fig:L1}
for the configuration of the curves
$C_{i,j}, \Gamma_i$ and $\Delta_i$.
We note that as basis of the cohomology group
$H^2(E_i,\ZZ)$ we can take the following curves:
\begin{itemize}
\item $C_{1,1}, \,\Delta_1,\,\Gamma_1$
and $\ol\Gamma_{n-1}$
when $i=1$,
\item $C_{i,i}, \,\Delta_i,$ and $\Gamma_i$
when $1<i<n-1$,
\item $C_{n-1,n-1}$ and $\Gamma_{n-1}$ 
when $i=n-1$.
\end{itemize}
In particular the restriction of 
the line bundle $\mathscr L_1$ to 
the divisor $E_i$ can be detected from
the intersection numbers with these curves. 
\begin{lemma}\label{lemma:int1}
The intersection numbers of $\mathscr L_1$ with the
above curves are given by
\begin{equation}\label{int01}
(\mathscr L_1, C_{i,i})_{Z_1} = 
\left\{
\begin{array}{cr}
0, &  i= 1,2,n-1,\\
-1, & 2<i<n-1,
\end{array}
\right.
\end{equation}
\begin{equation}\label{int02}
(\mathscr L_1, \Delta_i)_{Z_1} = 
\left\{
\begin{array}{cl}
0, &  i= 1,\\
1, & 1<i<n-1,\\
n-3, & i=n-1,
\end{array}
\right.
\end{equation}
and
\begin{equation}\label{int03}
(\mathscr L_1, \Gamma_i)_{Z_1} = 
\left\{
\begin{array}{cl}
n-2-i, &  1\le i\le n-2,\\
0, & i=n-1.
\end{array}
\right.
\end{equation}
\end{lemma}

\proof
First noting the relation $\mu^*F \sim 
f_1^*\mathscr O_{\Lmd}(1) + E$,
we have the following useful formula:
\begin{align}\label{L1-1}
\mathscr L_1 \sim f_1^*\mathscr O_{\Lmd}(n-2) + 
(E_1 + \ol E_1) + \sum_{j=2}^{n-1} (j-1)(E_j + \ol E_j).
\end{align}
Since the curves $C_{i,i}$ and $\Delta_i$ are contained in a fiber $S_i^+\cup S_i^-$
of $f_1$ and the curve $\Gamma_i$ is a section of $f_1$, we have
\begin{align}\label{int100}
(f_1^*\mathscr O(1), C_{i,i})_{Z_1} = 
(f_1^*\mathscr O(1), \Delta_i )_{Z_1} = 0,
\quad
(f_1^*\mathscr O(1),\Gamma_i)_{Z_1}
= 1.
\end{align}
For the intersection numbers of $E_i$ with 
the curves in the lemma,
if the curve is not contained in $E_i$ but intersects $E_i$,
the intersection number is one as they intersect transversally at a point.
On the other hand, for a curve which is contained in $E_i$ such as the curve $\Delta_i$
with $1<i<n-1$, noting $\Delta_i = S_i^+\cap E_i$
and $\Delta_i$ is contained in $S_i^+$ as a
$(-1)$-curve, we have
\begin{align}\label{int101}
(E_i,\Delta_i)_{Z_1}=
(\Delta_i,\Delta_i)_{S_i^+}=-1
\quad (1<i<n-1).
\end{align}
Similarly, noting that 
$C_{i,i}$ and $\Gamma_i$ are contained in $E_i$ as $(-1)$-curves for these $i$, we have
\begin{align}\label{int102}
(E_i,C_{i,i})_{Z_1}=
(E_i,\Gamma_i)_{Z_1}=-1, 
\quad 1<i<n-1.
\end{align}
Also, for the curves in the remaining 
components $E_1$ and $E_{n-1}$, we have
\begin{align}\label{int103}
(E_1,C_{1,1})_{Z_1} = 1-n,\quad
(E_1,\Delta_1)_{Z_1} =-1,\quad
(E_1,\Gamma_1)_{Z_1}= 
(E_1,\ol\Gamma_{n-1})_{Z_1}= 0,
\end{align}
\begin{align}
\label{int104}
(E_{n-1},C_{n-1,n-1})_{Z_1}=
(E_i,\Gamma_{n-1})_{Z_1}=-1.
\end{align}
We note that \eqref{int101}--\eqref{int104}
uniquely specify the normal bundle
$\mathscr O_{Z_1}(E_i)|_{E_i}$
for any $1\le i\le n-1$.
(These will be frequently used later.)

With these preparatory data,
the intersection numbers in the lemma can be computed readily.
For example, for proving \eqref{int01},
looking the right picture in Figure \ref{fig:L1},
when $i$ satisfies $1<i<n-1$, 
the curve $C_{i,i}$ intersects only with $E_{i-1}$ and $E_i$.
Therefore using \eqref{L1-1} with the aid of
\eqref{int100} and \eqref{int102}, when $2<i<n-1$, we compute
\begin{align}
(\mathscr L_1,C_{i,i})_{Z_1} & = 
(i-2) ( E_{i-1}, C_{i,i})_{Z_1} +
(i-1) ( E_{i}, C_{i,i})_{Z_1} \notag\\
&= (i-2)+(i-1)(-1)=-1.\notag
\end{align}
(The case $i=2$ requires an independent treatment
because of the form of the R.H.S.\,of \eqref{L1-1}.)
Hence we obtain the second case in \eqref{int01}.
The first case in \eqref{int01} can be obtained in 
a similar way by using \eqref{int102}--\eqref{int104}.
The other two assertions \eqref{int02} and \eqref{int03}
can be obtained in a similar way by using \eqref{L1-1}--\eqref{int104}.
\proofend
%
%
%
%

\begin{remark}\label{rmk:triv}
{\em
We have $(\mathscr L_1,C_{n-1,n-1})_{Z_1}
=(\mathscr L_1,\Gamma_{n-1})_{Z_1}=0$ by the lemma. 
It follows that $\mathscr L_1$ is trivial over
the two components $E_{n-1}$ and $\ol E_{n-1}$,
which will turn out to be useful later.
This is a reason why we choose 
the particular small resolutions $Z_1\to\hat Z$.
}
\end{remark}
By using the lemma, we show the following 
proposition which will be needed for proving Proposition \ref{prop:Z}
(iii). 
\begin{proposition}\label{prop:restsect}
For the restriction of the line bundle $\mathscr L_1$ over $Z_1$,
we have
$$h^0 \big((S_1\sqcup\cdots \sqcup S_{n-2})\cup E,\mathscr L_1\big) 
= n.$$
\end{proposition}

\proof
We first show $h^0(E, \mathscr L_1)=2$.
Recalling that $\mathscr L_1$ is trivial 
over $E_{n-1}$ and $\ol E_{n-1}$,
we prove that the restriction map
\begin{align}\label{rest5}
H^0(E,\mathscr L_1) \,\lra\, H^0(E_{n-1}\sqcup\ol E_{n-1},\mathscr L_1)\simeq\mathbb C^2
\end{align}
 is isomorphic, by verifying that 
any section over $E_{n-1}\sqcup \ol E_{n-1}$  extends in a unique way to the whole $E$.
For this we consider  
the following three restriction maps:
\begin{enumerate}
\item $H^0(E_i,\mathscr L_1) \lra H^0 (\Gamma_i,\mathscr L_1),
\,3\le i\le n-2$,
\item
$H^0(E_1,\mathscr L_1)\lra H^0(\ol{\Gamma}_{n-1},\mathscr L_1)$,
\item
$H^0(E_2,\mathscr L_1) \lra H^0(\Gamma_1\sqcup \Gamma_2,
\mathscr L_1).$
\end{enumerate}
By using 
Lemma \ref{lemma:int1}, it is elementary
to see that all of these maps are {\em isomorphic}.
Then  by the isomorphicity of the  
restriction maps (1) and (2) and their conjugations by the real structure,
 any section $s\in H^0(E_{n-1},\mathscr L_1)$
 (resp.\,$t\in H^0(\ol E_{n-1},\mathscr L_1)$) successively extends in a unique way to a section 
over the connected union $E_3\cup E_4\cup \cdots
\cup E_{n-1}\cup \ol E_1$
 (resp.\,$\ol E_3\cup \ol E_4\cup \cdots
\cup \ol E_{n-1}\cup E_1$).
 In particular, any section $(s,t)\in H^0(E_{n-1}\sqcup\ol E_{n-1},\mathscr L_1)$ uniquely determines a section over the curves
 $\Gamma_1\sqcup \Gamma_2$ and
  $\ol\Gamma_1\sqcup \ol\Gamma_2$.
Hence by the isomorphicity of (3), 
we conclude that $(s,t)\in H^0(E_{n-1}\sqcup\ol E_{n-1},\mathscr L_1)$ determines 
a section over the remaining components $E_2$ and $\ol E_2$ 
in a unique way.
Therefore the restriction map \eqref{rest5} is isomorphic,
and we obtain $h^0(E,\mathscr L_1)=2$.

On the other hand, 
recalling that for any $1\le k\le n-2$ the surface $S_k$ in $Z_1$ is canonically isomorphic to 
the original surface $S_k$ in $Z$,
from the definition of the line bundle $\mathscr L_1$,
the restriction  $\mathscr L_1|_{S_k}$ is linearly
equivalent to the movable part of
the  system
$|(n-2)K_{S_k}^{-1}|$.
Therefore we have 
$h^0(S_k,\mathscr L_1)= 3$ by Proposition 
\ref{prop:S} (ii).
We now show that the restriction map
\begin{align}\label{rest6}
H^0(S_k,\mathscr L_1) \,\lra\, H^0(S_k\cap E,\mathscr L_1)
\end{align}
is surjective for any $1\le k\le n-2$.
For this we note that under the above isomorphism, 
the intersection $S_k\cap E$ is  identified with the 
cycle $C$ on $Z$.
Therefore from the self-intersection numbers \eqref{B1} of 
the components of the cycle $C$,
we can readily see that the degree of the line bundle $\mathscr L_1|_{S_k}$
over each component of $S_k\cap E$ satisfies
\begin{equation}\label{int04}
(\mathscr L_1, S_k\cap E_i)_{Z_1} = 
\left\{
\begin{array}{cl}
0, &  i\neq 2,\\
1, & i=2.
\end{array}
\right.
\end{equation}
From these it is elementary to show that 
$
h^0(S_k\cap E,\mathscr L_1) = 2$ and $h^0(S_k, \mathscr L_1|_{S_k}-(S_k\cap E)) = 1$.
As the latter space is exactly the kernel of the restriction map \eqref{rest6},
by  dimension counting, we obtain  that \eqref{rest6} is surjective.

For completing the proof of the proposition, let 
\begin{align}\label{discrep1}
 H^0(E,\mathscr L_1) \oplus
\bigg(
\bigoplus_{k=1}^{ n-2}H^0(S_k,\mathscr L_1)
\bigg)
\,\stackrel{d}{\lra}\,
\bigoplus_{k=1}^{ n-2}H^0(E\cap S_k,\mathscr L_1)
\end{align}
be the linear map which takes differences on all connected components of the intersection
$E\cap (S_1\sqcup \cdots \sqcup S_{n-1})$.
Then since the  divisor \eqref{restdiv} is smooth normal crossing, we have
$$
\Ker\, d \simeq H^0 \big((S_1\sqcup\cdots \sqcup S_{n-2})\cup E,\mathscr L_1\big).
$$ 
The map $d$ is surjective since all the maps \eqref{rest6} are surjective.
Therefore, again by dimension counting, we finally obtain
\begin{align*}
h^0 \big((S_1\sqcup\cdots \sqcup S_{n-2})\cup E,\mathscr L_1\big)
&= 2 + 3(n-2) - 2(n-2)
= n,
\end{align*}
and we obtain Proposition \ref{prop:restsect}.
\proofend

\vsp
Continuing a proof of Proposition \ref{prop:Z} (iii),
we define another line bundle over $Z_1$ by
\begin{align}\label{L1'}
\mathscr L_1' := \mathscr L_1 - \sum_{k=1}^{n-2} S_k - E.
\end{align}
This is the kernel from the restriction of $\mathscr L_1$ to the key divisor \eqref{restdiv}.
This line bundle is still real.
For this line bundle we have the following critical vanishing result:
\begin{proposition}\label{prop:mainvan}
Let $\mathscr L_1'$ be the line bundle
\eqref{L1'} over $Z_1$ as above. Then we have
$$h^0(Z_1,\mathscr L_1') = 1\,\text{ and }\, H^q(Z_1,\mathscr L_1') = 0 
\,\text{ for any\, $q>  0$}.$$
\end{proposition}

\proof
First, summing up the relation $S_k  \sim \mu^*F - E$  on $Z_1$ for each $1\le k\le n-2$, we have $\sum_{k=1}^{n-2} S_k  \sim \mu^*((n-2)F) - (n-2)E$.
Hence from \eqref{L} we readily obtain
\begin{align}\label{kernel1}
\mathscr L_1'  \sim \sum_{i=3}^{{n-1}} 
(i-2)(E_i + \ol E_i).
\end{align}
We note that thanks to Lemma \ref{lemma:int1} and the intersection numbers \eqref{int101}--\eqref{int104}, all the restrictions $\mathscr O_{Z_1}(E_i + \ol E_i)|_{E_j\sqcup \ol E_j}$ are explicitly computable for any $i$ and $j$. 
By using these, it is possible to decrease the coefficients
in \eqref{kernel1} one by one  without changing 
arbitrary cohomology groups,
by subsequently considering the restrictions to the exceptional divisors 
in the following order:
$$
\begin{matrix}
E_3\sqcup \ol E_3 & E_4 \sqcup \ol E_4 & E_5 \sqcup \ol E_5 & 
E_6\sqcup \ol E_6\\
    &     &     & 1\\
    &     &   2 & 3\\
    &  4  &   5 & 6\\
 7  &  8  &   9 & 10\\  
\end{matrix}
$$
Here for simplicity we are displaying the order in the case $n=7$.
This means that we first restrict $\mathscr L_1'$ to $E_{6}\cup
\ol E_{6}$,
and second restrict $\mathscr L_1' - (E_{6}+\ol E_{6})$ 
to $E_{5}\cup\ol E_{5}$, and third
restrict $\mathscr L_1' - (E_{6}+\ol E_{6}+E_{5} + \ol E_{5})$ to $E_{4}\cup\ol E_{4}$, and so on.
Then if we  write $\mathscr O_{E_i}(a,b)$
for the pullback of the line bundle
$\mathscr O(a,b)$ over the original $E_i\simeq\CP^1\times\CP^1\subset \hat Z$ by the small resolution $Z_1\to \hat Z$,
then the  bundles over $E_i$ appearing
in the above restriction process are always of the form
$\mathscr O_{E_i}(-1,d)$ for some $d\in\mathbb Z$.
(Here, more precisely, $\mathscr O(0,1)$ represents
the fiber class of the projection $E_i\to \Lmd$.)
By the real structure we have the same result for the restrictions to $\ol E_i$.
Therefore any cohomology group vanishes, and we finally obtain 
$$
H^q (Z_1,\mathscr L_1')\simeq H^q(Z_1,\mathscr O_{Z_1}),
\quad q\ge 0.
$$
Further since the morphism $\mu:Z_1\to Z$ is birational, we have $H^q(\mathscr O_{Z_1})
\simeq H^q(\mathscr O_Z)$ for any $q\ge 0$.
Therefore we have $H^q(\mathscr L_1')\simeq H^q(\mathscr O_Z)$ for 
any $q\ge 0$, which easily implies the claim of the proposition.
\proofend

\vsp
\noindent
{\em Proof of Proposition \ref{prop:Z} (iii).}
From the definition of the line bundles $\mathscr L_1$ and $\mathscr L_1'$, we have 
an isomorphism 
$H^0(Z, (n-2)F) \simeq H^ 0 (Z_1, \mathscr L_1)
$
and an exact sequence
\begin{align}\label{ses2}
0 \,\lra\, \mathscr L_1' \,\lra\, \mathscr L_1 \,\lra\, \mathscr L_1 |_{(S_1\sqcup\cdots \sqcup S_{n-2})\,\cup \,E} \,\lra\, 0.
\end{align}
Hence by making use of Propositions
\ref{prop:restsect} and \ref{prop:mainvan}, we obtain 
$$
h^0(Z, (n-2)F) = h^0 (\mathscr L_1 |_{(S_1\sqcup\cdots \sqcup S_{n-2})\,\cup \,E}) + 1 = n+1.
$$
This completes a proof.
\proofend

\vsp
This readily implies the following.
\begin{proposition}\label{prop:surj}
For any non-singular member $S\in |F|$,
the restriction map 
$H^0(Z,(n-2)F) \to H^0((n-2)K_S^{-1})$
(see \eqref{rest1}) is surjective.
\end{proposition}

\proof
From the standard exact sequence
$$
0 \,\lra\, (n-3)F \,\lra\, (n-2)F \,\lra\, 
(n-2)K_S^{-1} \,\lra\, 0
$$
and Proposition \ref{prop:Z} (ii), we obtain an exact sequence
\begin{align}\label{les1}
0 \,\lra\, \CC^{n-2} \,\lra\, H^0((n-2)F) \,\lra\, 
H^0\big((n-2)K_S^{-1}\big).
\end{align}
Moreover we have $h^0((n-2)F) = n+1$
by Proposition \ref{prop:Z} (iii), and 
also $h^0 ((n-2)K_S^{-1}) = 3$ by Proposition \ref{prop:S} (ii).
Therefore from the exact sequence \eqref{les1} we obtain that
the restriction map is surjective.
\proofend

\vsp
For the rational map associated to the system $|(n-2)F|$, 
we have the following

\begin{proposition}\label{prop:Z2}
Let $\Phi:Z\to \CP^n$ be the rational map associated
to the linear system $|(n-2)F|$.
Then we have: 
\begin{enumerate}
\item[(i)] 
The image $\Phi(Z)$ is a scroll of 2-planes over
a rational normal curve in $\CP^{n-2}$.
\item[(ii)] The map $\Phi$ is two to one over the scroll.
\item[(iii)] The branch divisor of $\Phi$ is a cut of 
the scroll by a single quartic hypersurface.
\end{enumerate}
\end{proposition}

\proof
We just write an outline of the proof,
since these can be proved in a similar way to
\cite[Propositions 3.2 and 3.4]{HonDS4_1}.
From the subspace $S^{n-2}H^0(F)\subset H^0((n-2)F)$
and the rational maps associated to these linear systems,
we obtain the following commutative diagram of rational maps:
\begin{equation}\label{016}
 \CD
 Z@>\Phi >> \CP^n\\
 @V\Phi_{|F|} VV @VV\pi V\\
 \CP^1@>{\iota}>> \CP^{n-2},
 \endCD
 \end{equation}
where $\pi$ is the  linear projection induced from the above 
inclusion of the subspace, and $\iota$ is an embedding
as a rational normal curve.
Moreover by Proposition \ref{prop:surj}, 
the restriction map \eqref{rest1} is  surjective.
Hence the restriction of $\Phi$ to any non-singular member $S\in|F|$ 
is exactly the rational map $\phi:S\to \CP^2$ associated
to the net $|(n-2)K_S^{-1}|$.
Hence by the commutativity of the diagram \eqref{016} we obtain
the claim
(i).
Also, Proposition \ref{prop:S} (v) means 
the assertions (ii) and (iii).
\proofend


\section{Finding reducible members}\label{s:red}
Our final goal is to determine a defining equation of the 
quartic hypersurface which cuts out the branch divisor of
the map $\Phi:Z\to Y$ (see Proposition \ref{prop:Z2} (iii)).
For this purpose, in this section, we find two reducible members
of the system $|(n-2)F|$, each of which consists of two irreducible components.
As in the  case of $4\CP^2$ studied in \cite{HonDS4_1,HonDS4_2},
existence of these reducible members brings a strong constraint for a defining equation of the quartic hypersurface.
But in contrast with the case of $4\CP^2$,
for many reasons,
finding these divisors in the present case is incomparably difficult.

Let $S$ be the rational surface constructed in Section \ref{ss:S},
which is contained in the twistor space $Z$ 
as a real member of $|F|$ by our assumption.
Let $\epsilon:S\to\CP^1\times\CP^1$ be the composition of the explicit blowups
give in Section \ref{ss:S}.
(So $\epsilon$ is the composition of 
$S\to S_0$ and $S_0\to \CP^1\times\CP^1$.)
Let $e_1,\cdots,e_n$ and $\ol e_1,\cdots,\ol e_n$ be
the elements of $H^2(S,\mathbb Z)$ which are represented by
the exceptional curves of $\epsilon$,
named after the following natural rule:
$e_1$ and $e_2$ are represented by the exceptional curves
of the two blownup points on $C_1$ for obtaining
the surface $S_0$,
and $e_3$ is represented by the exceptional curve
of the blownup point on $C_2$
 for obtaining $S_0$.
The remaining classes $e_4, e_5 \cdots, e_n$ are chosen
in a standard way from the iterated blowup $S\to S_0$.
In particular, the classes $e_4,e_5,\cdots, e_{n-1}$ are not represented by
an irreducible curve, and $e_n$ is exactly the class of the curve $C_{n-1}$.
From the choice these classes satisfy $(e_i, e_j)_S=-\delta_{ij}$.
As a basis of the cohomology group 
$H^2(S,\ZZ)$ we can take the following $(2n+2)$ classes:
\begin{align}
\label{basis01}
e_1,\cdots,e_n,\ol e_1,\cdots,\ol e_n,
\epsilon^*\mathscr O(1,0),\,
 \epsilon^*\mathscr O(0,1).
\end{align}
Note that the roles of the first two classes $e_1$ and
$e_2$ are in some sense 
`symmetric'.

Let $\varpi:Z\to n\CP^2$ be the twistor fibration, and
let $\alpha_1,\cdots,\alpha_n$ be elements of $H^2(n\CP^2,\mathbb Z)$ which are uniquely determined from the condition 
$(\varpi^*\alpha_i)|_S = e_i - \ol e_i$ (see \cite{PP94}
for the structure of the restriction $\varpi |_S : S\to n\CP^2$).
As the pullback $\varpi^*:H^2(n\CP^2,\mathbb Z)\to
H^2(Z,\mathbb Z)$ is injective, in the following we just write $\aaa_i$ to mean 
$\varpi^*\aaa_i$.
Then the purpose of this section is to prove the following existence result:

\begin{proposition}
\label{prop:red}
Let $\mathscr M$ be the holomorphic line bundle over $Z$ whose cohomology class is given by
\begin{align}\label{m}
\frac{n-2}2 F - \frac12 \bigg\{
(n-2) \aaa_1 + (n-4) \sum_{i=2}^{n} \aaa_i
\bigg\}.
\end{align}
Then the linear system $|\mathscr M|$ consists of a single member,
and it is irreducible.
Also, the same conclusion holds for another cohomology class
\begin{align}\label{m'}
\frac{n-2}2 F - \frac12 \bigg\{
(n-2) \aaa_2+ (n-4) \sum_{i\neq 2} \aaa_i
\bigg\}.
\end{align}
\end{proposition}
We note that the latter class \eqref{m'} is obtained from \eqref{m} by
just exchanging the role of $\aaa_1$ and $\aaa_2$.
(This reflects the above `symmetric' property of $e_1$ and $e_2$.) 
We also note that the line bundle $\mathscr M$ and
the other one are {\em not} real,
and satisfy the relation 
$$
\mathscr M +\ol{ \sigma^*\mathscr M } \,\simeq (n-2)F.
$$
Therefore the single member of $|\mathscr M|$ (and also the single member of another system)
gives a reducible member of the system $|(n-2)F|$
consisting  of two irreducible components.

Our proof of Proposition \ref{prop:red} broadly proceeds in a similar way to
Proposition \ref{prop:Z} (iii).
Namely we pullback the line bundle $\mathscr M$ to the 
same blownup space $Z_1$, subtract obvious fixed components
from the pullback, and then
restrict the resulting bundle to some divisors of smooth normal crossing.
But the choice of the last divisor is much more subtle than
we did in Section \ref{s:ps}
as we  see below.

We begin with determining fixed components of 
the linear system $|\mathscr M|_S|$ on the surface $S$:

\begin{proposition}
\label{prop:fc}
Let $S\in |F|$ be any non-singular member of the pencil $|F|$.
Then the linear system $|\mathscr M|_S|$ contains the following curve as  fixed
components at least:
\begin{align}\label{bc2}
(n-3)C_1 + \sum_{i=2}^{ n-2} (n-1-i) C_i.
\end{align}
In other words, any section of the line bundle $\mathscr M|_S$ vanishes
along the curve $C_i$  by the order indicated by the coefficient at least.
\end{proposition}

\proof
Though this is not immediate to see,
it can be proved in an elementary way,
so we just give an outline.
From the explicit form of the line bundle $\mathscr M$ and 
the relation $\aaa_i|_S= e_i - \ol e_i$ for $1\le i\le n$,
we can concretely write down the cohomology class of 
the line bundle $\mathscr M|_S$, in terms of the basis \eqref{basis01}.
Also, the cohomology class of the curve $C_i$ can be 
expressed in terms of the same basis.
Therefore we can compute the intersection numbers of the line bundle
$\mathscr M|_S$ with the curve $C_i$.
From this, by checking negativity
or vanishing of the intersection numbers successively,
we can show that any section of $\mathscr M|_S$ 
has to vanish along the curve \eqref{bc2} with multiplicities indicated by the coefficients. 
\proofend


\vsp
Let $\mu:Z_1\to Z$ be the birational morphism given in Section \ref{ss:ps}.
By Proposition \ref{prop:fc}, if we define a line bundle
$\mathscr M_1$ over $Z_1$ by
\begin{align}\label{m1}
\mathscr M_1:= \mu^* \mathscr M - 
(n-3)E_1 - \sum_{i=2} ^{ n-2} (n-1-i) E_i,
\end{align}
then we have an isomorphism
\begin{align}\label{iso77}
H^0(Z,\mathscr M)\simeq
H^0(Z_1,\mathscr M_1).
\end{align}
We are going to compute the right-hand-side by restricting $\mathscr M_1$ to 
the divisor 
\begin{align}\label{restdiv2}
E + \sum_{i=1}^{n-2} S_i^-,
\end{align}
where as in Section \ref{ss:ps}, $S_i^-$ is the strict transform of an irreducible component of a reducible
member of the pencil $|F|$, and $E$ is the total sum
of the exceptional divisors of the birational morphism $\mu$.
We note that the divisor \eqref{restdiv2} is again smooth normal crossing, and that 
the degree of the divisor \eqref{restdiv2} is equal to
that of $\mathscr M$.
(Note that in \eqref{restdiv2} the divisor $S_{n-1}^-$ is {\em not} included.
 The reason why we restrict to this particular divisor among 
numerous possible choices
would become evident in the course of the proof of 
Proposition \ref{prop:mainvan3} below.)

We define another line bundle $\mathscr M'_1$ over $Z_1$ by
\begin{align}
\mathscr M'_1 := \mathscr M_1 - \bigg(E +  \sum_{i=1} ^{n-2} S_i^-\bigg),
\end{align}
which is  the kernel of the restriction
of the line bundle $\mathscr M_1$ to the divisor
\eqref{restdiv2}.
Then
one of the keys for computing $H^0(\mathscr M_1)$ is the following 
critical vanishing result:
\begin{proposition}\label{prop:mainvan3}
For any $q\ge 0$, we have\,
$
H^q(Z_1,\mathscr M_1') = 0.
$
\end{proposition}  

\proof
The strategy is the same as a similar result Proposition \ref{prop:mainvan} in the last section,
but the required computations are more involved.
In this proof for distinguish 
divisor $S_i^-$ in $Z$ and its strict transform into $Z_1$,
we write $\mathscr O_{Z}(S_i^-)$ and $\mathscr O_{Z_1}(S_i^-)$ respectively.

First by using the fact that the restriction map
$H^2(Z,\mathbb Z)\to H^2(S,\mathbb Z)$ is injective,
and also from a concrete form of the divisor $S_i^-|_S$
which is exactly a half of the cycle $C$,
 it is possible to write down the Chern classes of the original divisors
$S_i^-$ in $Z$, in terms of $F$ and the classes $\aaa_1,\cdots,\aaa_n$.
The result is as follows:
\begin{align}\label{chern1}
\mathscr O_Z(S_i^-) \simeq \frac 12 F - \frac 12
\sum _{j=1}^{n}\eee_{ij}\aaa_j, \,\,{\text{ for }}1\le i\le n-2,
\end{align}
where $\eee_{ij}=1$ except $j= n-i+1$ while
$\eee_{ij}=-1$ if $j= n-i+1$, and for $i=n-1$
\begin{align}\label{chern2}
\mathscr O_Z(S_{n-1}^-) \simeq \frac 12 F - \frac 12
\sum _{j=1}^n \aaa_j.
\end{align}
(The formula \eqref{chern2} will be needed later.)
Summing up \eqref{chern1} for $1\le i\le n-2$,
we easily obtain 
\begin{align}\label{sum1}
\sum_{i=1}^{n-2} \mathscr O_Z(S_i^-) \simeq
\frac{n-2}{2}F - \frac 12\bigg\{
(n-2)(\aaa_1+\aaa_2) + (n-4) \sum_{i=3}^{n}
\aaa_i\bigg\}.
\end{align}

Next as the divisor $S_1^-$ in $Z$ contains
the $n$ curves $\ol C_2,\ol C_3,\cdots, \ol C_{i-1}$ and $C_1$ by  multiplicity one,
we have
$\mathscr O_{Z_1}(S_1^-) \simeq \mu^*\mathscr O_Z(S_1^-)
- \sum_{i=2}^{ n-1}\ol E_i - E_1$. 
We have a similar isomorphism for the line bundle
$\mathscr O_{Z_1}(S_i^-)$ for any $i$.
Summing these up for $1\le i\le n-2$, we obtain
\begin{align}\label{sum2}
\sum_{i=1}^{n-2} \mathscr O_{Z_1}(S_i^-)
\simeq 
\mu^* \bigg(\sum_{i=1}^{n-2} \mathscr O_Z (S_i^-)\bigg)
- \bigg(
\sum_{i=1}^{n-1} (n-1-i) E_i +
\sum_{i=1}^{n-1}  (i-1) \ol E_i
\bigg).
\end{align}
After substituting \eqref{sum1} into \eqref{sum2}, we deduce
\begin{align}
\notag\mathscr M_1' & = 
\mu^* \mathscr M - 
(n-3)E_1 - \sum_{i=2}^{ n-2} (n-1-i) E_i 
-E - \big({\text{R.H.S. of \eqref{sum2}}}\big)\\
&= \mu^* \mathscr O_Z(\aaa_2) - \sum_{i=2}^{ n-1}
E_i + \sum_{i=1}^{n-1}(i-2)\ol E_i.
\label{cancel}
\end{align}
 (In the equality \eqref{cancel} 
almost all terms in the pullback term canceled out, and this is the reason why we
choose the particular divisor \eqref{restdiv2} for the restriction.
The `smallness' of the pullback term is crucial as we  see in the following argument.)

So for the proof of the proposition it suffices to show that the cohomology group $H^q$ 
of the line bundle \eqref{cancel} vanishes for any $q$.
We first see that  the first summation
in \eqref{cancel} can be entirely removed without changing
any cohomology group.
By adding $E_2$, we obtain the standard exact sequence
\begin{align}\label{ses3}
0 \,\lra\, \mathscr M_1' \,\lra\,
\mathscr M'_1+ E_2 \,\lra\,
\mathscr M'_1+ E_2\,|_{E_2} \,\lra\,0.
\end{align}
For the restricted term, from \eqref{cancel} we have
$\mathscr M'_1+ E_2\,|_{E_2} \simeq 
\mu^* \mathscr O_Z(\aaa_2)  - E_3 |_{E_2}$.
Further since the curves $e_2\cup\ol e_2$ and $C_2$ (in $S$)
 are disjoint, this is isomorphic to $-E_3|_{E_2}$.
 It is easy to see that all cohomology groups vanish for the last class,
 and  hence we obtain isomorphisms $H^q(\mathscr M_1')
 \simeq H^q(\mathscr M_1'+E_2)$ for any $q\ge 0$.
 Repeating this process by  adding $E_3,E_4,\cdots, E_{n-1}$ one by one, 
 we finally obtain an isomorphism
\begin{align}\label{iso4}
H^q(\mathscr M_1')\simeq H^q
\bigg(
\mu^* \mathscr O_Z(\aaa_2) 
 + \sum_{i=1}^{n-1}(i-2)\ol E_i.
\bigg),
\quad q\ge 0.
\end{align}
(We remark that all the  line bundles over $E_i$-s appearing in this restriction process
are mutually isomorphic, which  considerably decreases the computations.)

Note that a negative term $-\ol E_1$ is still included in the ingredient of the 
R.H.S.\,of \eqref{iso4}.
We next show that this term can also be removed.
Recalling that in the surface $S$ the two curves
$\ol C_1$ and $e_2$ are disjoint and that $\ol C_1$ and $\ol e_2$ intersect transversally at a point,
we have
\begin{align}
\big(\mu^*\mathscr O_Z (\aaa_2)\big)|_{\ol E_1} 
 \simeq \mu^*\big(\mathscr O_Z(\aaa_2)|_{\ol C_1}\big)
&\simeq \mu^* \big(
\mathscr O_S(\aaa_2)|_{\ol C_1}\big)\notag\\
&\simeq \mu^*\big(\mathscr O_S(e_2-\ol e_2)|_{\ol C_1}\big)
\simeq 
\mu^*\big(\mathscr O_{\ol C_1}(-1)\big).
\end{align}
Hence noting that the component $\ol E_3$ 
is not included in the R.H.S.\,of \eqref{iso4},
we obtain that the restriction of [the ingredient of the R.H.S.\,of \eqref{iso4}
plus $\ol E_1$]
to  the divisor $\ol E_1$ is  isomorphic to just
$
\mu^*\big(\mathscr O_{\ol C_1}(-1)\big),
$
whose all cohomology groups can be easily seen to vanish.
Hence by the exact sequence similar to \eqref{ses3},
we can remove the negative term $-\ol E_1$
 without changing any cohomology group.

Thus for completing the proof
of Proposition \ref{prop:mainvan3} we are reduced to show 
\begin{align}\label{van3}
H^q \bigg(\mu^* \mathscr O_Z(\aaa_2) + \sum_{i=3}^{n-1}(i-2)\ol E_i \bigg) = 0,
\quad q\ge 0.
\end{align}
Here we note that  the summation  in \eqref{van3} is 
exactly the one included in the line bundle $\mathscr L_1'$
(see \eqref{kernel1}) in the last section.
Therefore the computations in the proof of 
Proposition \ref{prop:mainvan} perfectly work in order to decrease
the coefficients of $\ol E_i$-s one by one,
and finally we obtain an isomorphism
\begin{align}\label{iso3}
H^q \bigg(\mu^* \mathscr O_Z(\aaa_2) + \sum_{i=3}^{ n-1}(i-2)\ol E_i \bigg)\simeq
H^q \big(\mu^* \mathscr O_Z(\aaa_2)\big),
\quad q\ge 0.
\end{align}
The R.H.S. of \eqref{iso3} is of course isomorphic to 
$H^q(Z,\mathscr O_Z(\aaa_2))$.
For $q=0$ and $q=3$, this is zero by obvious reasons.
For $q=2$, this also vanishes by the vanishing theorem 
of Hitchin \cite{Hi80}.
On the other hand, the Riemann-Roch formula gives
\begin{align}\label{RR}
\chi (\mathscr O_Z(\aaa_2)) = 
\frac 16 \aaa_2^3 + \frac 14 \aaa_2^2  c_1
+ \frac1{12} \aaa_2 (c_1^2 + c_2) + \frac 1{24} c_1c_2,
\end{align}
where $c_i$ denotes the Chern class of $Z$.
We have $\aaa_2^3 = 0$ since $\aaa_2$ is a lift from $n\CP^2$.
We have $\aaa_2^2\cdot c_1 = -4$ because $K_Z$ is of degree 4
over a twistor line.
On the other hand both $c_1^2$ and $c_2$ are lifts from 
$n\CP^2$ (see \cite{Hi81}), and therefore their product with $\aaa_2$ is zero.
Finally $c_1c_2$ is  $24$.
Hence we obtain $\chi(\mathscr O_Z(\aaa_2)) = -1 + 1 = 0.$
Thus we get $H^1 (\mathscr O_Z(\aaa_2)) = 0$.
Therefore we obtain $H^q(\mathscr O_Z(\aaa_2)) = 0$ for any $q\ge 0$,
and finally obtain 
$H^q (\mathscr M'_1) = 0$ for any $q\ge 0$.
\proofend

\vsp

The following result is also indispensable for 
proving Proposition \ref{prop:red}.
\begin{proposition}\label{prop:nonvan2}
For the restriction of the line bundle $\mathscr M_1$ over $Z_1$, we have
$$h^0\big((S_1^-\sqcup\cdots\sqcup S_{n-2}^-)\cup E,\mathscr M_1\big) = 1.$$
\end{proposition}

For the proof, we first show the following
\begin{proposition}\label{prop:nonvan3}
We have \,$h^0(E,\mathscr M_1) = 1$.
\end{proposition}

\proof
The idea is similar to the first half of the 
proof of Proposition \ref{prop:restsect},
but since the line bundle $\mathscr M$ possesses terms coming from $n\CP^2$, the computations are much more involved.
As in the case of the line bundle $\mathscr L_1$, 
we exhibit
the restrictions of $\mathscr M_1$ to
the components of $E$
in terms of the basis of $H^2(E_i,\ZZ)$ given
just before Lemma \ref{lemma:int1}.

First for obtaining the restrictions of 
the pulled-back term $\mu^*\mathscr M$,
we first compute the intersection numbers
$(\mathscr M,C_i)_Z$ and $(\mathscr M,\ol C_i)_Z$.
(Note that since $\mathscr M$ is non-real, these
are not necessarily equal.)
Putting 
$
\aaa := (n-2) \aaa_1 + (n-4) ( \aaa_2+\aaa_3+\cdots+\aaa_n)
$
so that $2\mathscr M= (n-2) F - \aaa$,
and
taking a non-singular member $S\in |F|$, we  have
\begin{align}
2(\mathscr M,C_i)_Z & = 2(\mathscr M|_S,C_i)_S\notag\\
& = \big( (n-2)K_S^{-1},C_i \big)_S - (\aaa|_S, C_i \big)_S,\label{int200}
\end{align}
and a similar equality for $(\mathscr M,\ol C_i)_{Z}$.
As $C_i$ is a rational curve, we have $(K_S^{-1},C_i)_S = (C_i, C_i)_S +2$.
For computing  another term $(\aaa|_S,C_i)_S$,
it suffices to compute $(\aaa_j|_S,C_i)_S$ for each 
$1\le j\le n$,
and this is equal to $(e_j-\ol e_j,C_i)_S$
by our definition of the class $\aaa_j$.
From the choice of the classes
$e_1,\cdots,e_n$ given
at the beginning of this section,
it is not difficult to deduce  the relations
\begin{align}
e_j = \sum_{k=n-j+2}^{ n-1}\ol C_k,\quad
3\le j\le n-1,
\end{align}
and similar relations for $\ol e_j$ and $C_k$.
From these, by using the self-intersection numbers
\eqref{B1}, 
we can quickly compute
the intersection numbers $(e_j,C_i)_S$ and $(\ol e_j,C_i)_S$ for any $1\le i\le n-1$ and $1\le j\le n$.
By using these and \eqref{int200},
after long but elementary computations, we
obtain that 
\begin{equation}\label{int2}
(\mathscr M, C_i)_Z = 
\left\{
\begin{array}{cc}
-(n-2)(n-3) & i=1,\\
0 &  1<i<n-1,\\
1 & i=n-1,
\end{array}
\right.
\end{equation}
and
\begin{equation}
\label{int10}
(\mathscr M, \ol C_i)_Z = 
\left\{
\begin{array}{cc}
0 & 1\le i<n-1,\\
n-3 & i=n-1.
\end{array}
\right.
\end{equation}

On the other hand, since $E_i$ (resp.\,$\ol E_i$) is an exceptional divisor
over the curve $C_i$ (resp.\,$\ol C_i$) in $Z$,
we have
\begin{align}\label{rests1}
(\mu^*\mathscr M)|_{E_i} \simeq
\mu^* (\mathscr M|_{C_i})
\quad\text{and}\quad
(\mu^*\mathscr M)|_{\ol E_i} \simeq
\mu^* (\mathscr M|_{\ol C_i}).
\end{align}
Further, recalling the concrete form of the restriction of the birational morphism
$\mu$ to the divisor $E_i$ (see Section \ref{ss:ps}),
we have, for any $d\in\ZZ$,
\begin{align}\label{restrc1}
(\mu^*\mathscr O_{C_i}(d))|_{C_{i,i}} \simeq
\mathscr O_{C_{i,i}}(d),\quad 1\le i\le n-1,
\end{align}
and
\begin{align}\label{restrc2}
(\mu^*\mathscr M)|_{\Delta_i} \simeq
\mathscr O_{\Delta_i},
\quad
(\mu^*\mathscr M)|_{\Gamma_i} \simeq
\mathscr O_{\Gamma_i}.
\end{align}
Hence $(\mu^*\mathscr M, C_{i,i})_{Z_1}$ and
$(\mu^*\mathscr M, \ol C_{i,i})_{Z_1}$ are
exactly given by the R.H.S-s of \eqref{int2} and
\eqref{int10} respectively, and
$(\mu^*\mathscr M, \Gamma_i)_{Z_1} = 
(\mu^*\mathscr M, \Delta_i)_{Z_1} = 0$ for any $i$.

On the other hand, the intersection numbers of the subtraction terms 
(see \eqref{m1})
with the above curves in $E_i$ and $\ol E_i$
can be readily computed by using the intersection numbers \eqref{int101}--\eqref{int104}.

Combining these, by elementary calculations,
we can deduce that the intersection numbers
of $\mathscr M_1$ with the above curves are given by 
\begin{equation}
\label{int300}
(\mathscr M_1, C_{i,i})_{Z_1} = 
\left\{
\begin{array}{cc}
0 & i\in\{1,2,n-1\},\\
-1 & i\not\in\{1,2,n-1\},
\end{array}
\right.\quad
(\mathscr M_1, \ol C_{i,i})_{Z_1} = 
0\quad 1\le i\le n-1,
\end{equation}
\begin{equation}
\label{int301}
(\mathscr M_1, \Delta_i)_{Z_1} = 
\left\{
\begin{array}{cc}
0 & i\in\{1,n-1\},\\
1 & i\not\in\{1,n-1\},
\end{array}
\right.
(\mathscr M_1, \ol\Delta_{i})_{Z_1} = 
\left\{
\begin{array}{cc}
0 & i\neq n-1,\\
n-3 & i=n-1,
\end{array}
\right.\quad
\end{equation}
\begin{equation}
\label{int302}
(\mathscr M_1, \Gamma_i)_{Z_1} = 
\left\{
\begin{array}{cc}
n-i-2 & i\neq n-1,\\
0 & i=n-1,
\end{array}
\right.\quad
(\mathscr M_1, \ol\Gamma_i)_{Z_1} = 
0 \quad 1\le i\le n-1.
\end{equation}
In particular, we get that $\mathscr M_1$ is trivial over the $n$ components
$E_{n-1}, \ol E_1,\ol E_2,\cdots, \ol E_{n-1}$.
Therefore by connectedness (see Figure \ref{fig:L1}), we have 
$
h^0(E_{n-1}\cup  \ol E_1\cup \ol E_2\cup \cdots\cup  \ol E_{n-1},\mathscr M_1) = 1$.

For completing a proof of Proposition \ref{prop:nonvan3}, we show that 
the restriction map
\begin{align}\label{rest7}
H^0(E,\mathscr M_1) \,\lra\, 
H^0 \big(E_{n-1}\cup  \ol E_1\cup \ol E_2\cup \cdots\cup  \ol E_{n-1},\mathscr M_1\big )\simeq\CC
\end{align}
is isomorphic, by showing that any element
of the R.H.S.
uniquely extends to the whole $E$.
Similarly to the first part of the proof of Proposition \ref{prop:restsect},
we consider the following three restriction maps:
\begin{enumerate}
\item $H^0(E_i,\mathscr M_1) \lra H^0 (\Gamma_i,\mathscr M_1),
\quad 3\le i\le n-2$,
\item
$H^0(E_1,\mathscr M_1)\lra H^0(\ol\Gamma_{n-1},\mathscr M_1)$,
\item
$H^0(E_2,\mathscr M_1) \lra H^0(\Gamma_1 \sqcup \Gamma_2,
\mathscr M_1).$
\end{enumerate}
By using \eqref{int300}--\eqref{int302} it is elementary
to see that all these are isomorphisms.
Then by the argument in the proof
of Proposition \ref{prop:restsect},
we conclude that the restriction map \eqref{rest7} is isomorphic.
This means the claim of Proposition \ref{prop:nonvan3}.
\proofend

\vsp
Next for the proof of Proposition \ref{prop:nonvan2} we further need to show

\begin{proposition}\label{prop:nonvan4}
For the line bundle $\mathscr M$ over the  original twistor space $Z$, we have 
$$h^0(S_i^-,\mathscr M) = 1,\quad 1\le i\le n-2.$$
\end{proposition}

Before proceeding to the proof,
we note that on each divisor $S_i^-$ (resp.\,$S_i^+$)
there exist $(-1)$-curves $e'_1$ and $e'_2$ 
(resp.\,$\ol e'_1$ and $\ol e'_2$) 
such that $\aaa_j|_{S_i^- \cup S_i^+} = e'_j - \ol e'_j$
for $j=1,2$.
The existence of these $(-1)$-curves  can be derived from the self-intersection numbers of the components
of the cycle $C$ inside $S_i^-$ and $S_i^+$, 
and also from the fact that the intersection $S_i^-\cap S_i^+$ is a twistor line, which is contained in 
$S_i^+$ and $S_i^-$ as a $(+1)$-curve.
The curves $e'_j$ and $\ol e'_j$ are respectively homologous to 
the exceptional curves $e_j$ and $\ol e_j$ contained in 
each non-singular member $S\in |F|$.
Just like $e_1$ and $e_2$ in $S$, each of $e'_1$ and $e'_2$ (resp.\,$\ol e'_1$ and $\ol e'_2$) intersects $C_1$ (resp.\,$\ol C_1$) transversally
at a unique point respectively.

\vsp
\noindent
{\em Proof of Proposition \ref{prop:nonvan4}.}
Again we first compute the restriction of
the line bundle $\mathscr M$
to the divisor $S_i^-$ in a concrete form.
For this we need to compute the restriction of 
the class $\aaa_j$ to 
$S_i^-$ for any $j$, which is quite difficult in contrast with 
their restriction to $S$  (except the cases $j=1,2$).
To avoid this, we make use of the remaining component
$S_{n-1}^-$.
Let $\aaa$ be as in the proof of Proposition \ref{prop:nonvan3},
so that  $\mathscr M = \{(n-2)F - \aaa\}/2$.
By using the concrete form of $\aaa$ and the Chern class formula \eqref{chern2} of $S_{n-1}^-$,
we rewrite $\mathscr M$ as 
\begin{align}
\frac{n-2}2 F - \frac \aaa 2 &=
F + (n-4)  \bigg( \frac 12 F - \frac 12 \sum_{i=1}^{ n}
\aaa_i \bigg) - \aaa_1\notag\\
&= F + (n-4) S_{n-1}^- -\aaa_1.\label{tec}
\end{align}
Since $\aaa_1|_{S^-_i} = e_1'$ for any $i$
as above,
we obtain from \eqref{tec} that
\begin{align}\label{rest8}
\mathscr M|_{S_i^-} = F|_{S_i^-} + (n-4)S_{n-1}^-|_{S_i^-} - e'_1.
\end{align}
The first term $F|_{S_i^-}$ is exactly a half of the cycle $C$
contained in $S_i^-$, 
and the restriction $S_{n-1}^-|_{S_i^-}$
in the second term is also a part of the cycle $C$,
which can be immediately written down.
From these we obtain
\begin{align}
\label{rest9}
\mathscr M|_{S_i^-} = \sum_{j=i+1}^{n-1}\ol C_j
+ (n-3) \sum_{j=1}^i  C_j - e'_1.
\end{align}
By computing intersection numbers, it is immediate to see that 
the second term  $(n-3)\sum_{j=1}^i C_j$ is a  fixed component
of this system,
and that the system 
$\big|\sum_{j=i+1}^{n-1}\ol C_j\big|$ 
is a base point free pencil.
From the latter
we obtain that the system
$|\sum_{j=i+1}^{n-1}\ol C_j - e'_1|$
 has a
unique (effective) member, and that 
it  is disjoint from the cycle $C$.
These in particular mean the claim of the proposition.
\proofend

\vsp
By using Propositions \ref{prop:nonvan3} and \ref{prop:nonvan4}
 we show Proposition 
\ref{prop:nonvan2}:

\vspace{2mm}
\noindent
{\em Proof of Proposition \ref{prop:nonvan2}.}
First we compute the restriction of
the pulled-back bundle $\mu^*\mathscr M$ to the divisor
$S_i^-\subset Z_1$ for $1\le i\le n-2$.
For this we recall from Section \ref{ss:ps} that due to the small resolution $Z_1\to \hat Z$,
the restriction of the birational morphism
$\mu:Z_1\to Z$ to the divisor $S_i^-\subset Z_1$
is
not isomorphic but 
identified with the blowing-up at the point $\ol C_i\cap \ol C_{i+1}$,
and the curve $\ol\Delta_i$ is inserted as the
exceptional curve.
Then noting that, among the curves in the R.H.S. of \eqref{rest9},
only $C_{i+1}$ contains the blown-up point 
$\ol C_{i+1}\cap \ol C_i$,
and that the coefficient of $\ol C_{i+1}$ is one,
we have
\begin{align}\label{rest10}
(\mu^*\mathscr M)|_{S_i^-} = \sum_{j=i+1}^{n-1}\ol C_{i,j}
+ (n-3) \sum_{j=1}^i  C_{i,j} - e'_1 + \ol{\Delta}_i,
\end{align}
where we are using the curves $C_{i,j}$ defined in 
Section \ref{ss:ps}.
On the other hand, the subtraction term in $\mathscr M_1$ is
(see \eqref{m1})
\begin{align}\label{subtr}
(n-3)E_1 + \sum_{j=2}^{n-2} (n-1-j) E_j.
\end{align}
When $j>i$,
 the intersection $S_i^-\cap E_j$
is at most a point.
Hence by disposing them
we obtain that  the restriction of \eqref{subtr} to the divisor $S_i^-$ $(1\le i\le n-2)$ is given by
$$
(n-3)C_{i,1} + \sum_{j=2}^i (n-1-j) C_{i,j}.
$$
Subtracting this from \eqref{rest10} we get
\begin{align}\label{rest11}
\mathscr M_1|_{S^-_i} = \sum_{j=3}^i (j-2) C_{i,j} + 
\ol\Delta_i + \sum_{j=i+1}^{n-1} \ol C_{i,j}  - e'_1.
\end{align}
By computing intersection numbers, the first summation in 
\eqref{rest11} can easily seen to be fixed components
of $|\mathscr M_1|_{S_i^-}|$.
On the other hand from self-intersection numbers the system
$|\ol{\Delta}_i + \sum _{j=i+1} ^{n-1}\ol C_{i,j}|$ is again a base point free pencil,
and it follows that the system $|\ol{\Delta}_i + 
\sum _{j=i+1} ^{n-1} \ol C_{i,j} - e'_1|$ consists
of a single member.
In particular we obtain $h^0(S_i^-, \mathscr M_1) = 1$.
Moreover when $1\le i<n-1$, we have
$$
\bigg(\ol C_{i,n-1}, \sum_{j=i+1}^{n-1} \ol C_{i,j} 
+ \ol\Delta_i - e'_1
\bigg)_{S_i^-}
=
\bigg(\ol C_{i,n-1},\,  \ol C_{i,n-1} + \ol C_{i,n-2}  - e'_1
\bigg)_{S_i^-} 
= (-1)+ 1 - 0 = 0.
$$
When $i=n-1$, by replacing $\ol C_{i,n-2}$ with $\ol\Delta_{n-1}$,
we obtain the same conclusion.
Therefore
the unique member of the system $|\mathscr M_1|_{S_i^-}|$
is disjoint from the curve $\ol C_{i,n-1}$
for any $1\le i\le n-1$.

For completing the proof of Proposition \ref{prop:nonvan2}, 
by Proposition \ref{prop:nonvan3}, it is enough to show that 
any element of $H^0(E,\mathscr M_1)$ extends to 
$S_i^-$ for any $1\le i\le n-2$ in a unique way.
For this we first recall from the proof of Proposition \ref{prop:nonvan3}
that  $\mathscr M_1$ is trivial over the component $\ol E_{n-1}$,
and the restriction map 
$H^0(E,\mathscr M_1) \to H^0(\ol E_{n-1},\mathscr M_1)$ is isomorphic.
Further as above the unique member of the system $|\mathscr M_1|_{S_i^-}|$
is disjoint from the curve $\ol C_{i,n-1}$.
In particular we obtain that 
both of the two restriction maps
\begin{enumerate}
\item
$H^0(E,\mathscr M_1)
\,(\simeq\CC) \,\lra\, H^0(\ol C_{i,n-1},\mathscr M_1)\,(\simeq\CC)$,
\item
$H^0(S_i^-,\mathscr M_1)\,(\simeq\CC)  \,\lra\, H^0(\ol C_{i,n-1},\mathscr M_1)\,(\simeq\CC)$,
\end{enumerate}
are isomorphic.
As the divisor \eqref{restdiv2} is smooth normal crossing, this immediately means that 
any element of $H^0(E,\mathscr M_1)$ uniquely extends to
$S_i^-$ for any $1\le i\le n-2$.
\proofend

\vsp
Now we are able to prove the  proposition presented in the beginning of this section.

\vspace{2mm}
\noindent
{\em Proof of Proposition \ref{prop:red}.}
By the isomorphism \eqref{iso77},
in order to prove $h^0(Z,\mathscr M)=1$, it suffices to show $h^0(Z_1,\mathscr M_1) = 1$.
But now this is an immediate consequence of the standard exact sequence
\begin{align}\label{ses4}
0 \lra \mathscr M'_1 \lra \mathscr M_1 \lra \mathscr M_1 |_{(S^-_1\,\sqcup\cdots \sqcup\, S^-_{n-2})\,\cup E}\lra 0
\end{align}
and Propositions \ref{prop:mainvan3} and \ref{prop:nonvan2}.

Next let $D$ be  the unique member of 
the system $|\mathscr M|$ and  show that $D$ is irreducible.
Suppose that $D$ is reducible, and let $D_1$ be
any irreducible component of $D$.
Then we have $D_1+ \ol D_1\in |kF|$ for some $k$ with
$0<k<n-2$.
But by Proposition \ref{prop:Z} (i), (ii), we have
$|kF| = S^k H^0(F)$ for these $k$.
This means that $D_1$ is a degree-one divisor, or otherwise
$D_1\in |F|$.
Thus if $D$ is reducible, all irreducible components
must be some $S\in |F|$, $S_i^+$ or $S_i^-$.
In order to show that this cannot happen, 
we first notice that the coefficients
of $\aaa_1$ and $\aaa_2$ of 
the cohomology class $\mathscr M$ (see \eqref{m1}) do not
coincide (namely $(n-2)$ and $(n-4)$ respectively).
On the other hand  the class $F$ does not contribute
for the pullback term (i.e.\,$\aaa_i$-terms).
Furthermore, most importantly,
the Chern class formulae \eqref{chern1} and \eqref{chern2}
for degree-one divisors
imply that the coefficients of $\aaa_1$ and $\aaa_2$ 
coincide for any $S_i^+$ and $S_i^-$, $1\le i\le n-1$.
Therefore, 
$D$ cannot be a sum of 
$S_i^+$, $S_i^-$ $(1\le i\le n-1)$ and $S\in |F|$.

The claim for another line bundle \eqref{m'} follows
by exchanging the role of the two classes $\aaa_1$ and 
$\aaa_2$ in all the arguments throughout this section.
\proofend

\section{Elimination of the base locus of the pluri-half-anticanonical system}
\label{s:deg1}
In Section \ref{s:ps} we proved that the 
linear system $|(n-2)F|$ of the present twistor space induces 
a rational map $\Phi:Z\to Y\subset\CP^n$, where
$Y$ is the scroll of planes over the rational normal curve
$\Lambda\subset\CP^{n-2}$, and
that $\Phi$ is of degree two over $Y$
(Proposition \ref{prop:Z2}).
The restriction of the line bundle $(n-2)F$ to a smooth member $S\in |F|$ is 
isomorphic to $(n-2)K_S^{-1}$, and this line bundle has base points
along some components of the cycle $C$ (Proposition \ref{prop:S} (iii)).
Then from the surjectivity of the restriction map \eqref{rest1} (see Proposition \ref{prop:surj}),
we have a coincidence $\Bs\,|(n-2)F| = \Bs\,|(n-2)K_S^{-1}|$.
In this section, we give a complete elimination of this base locus,
via the space $Z_1$ we have used throughout Sections \ref{s:ps} and \ref{s:red}.
While the elimination requires some complicated calculations, this
process seems
to be indispensable for reaching our final goal.
We use the  notations from Sections \ref{s:ps} and \ref{s:red},
and continue to use the same letters to mean
subsets of $Z$ and their the strict transforms into $Z_1$.
Also, since all the operations preserve the real structure,
we often omit to mention the counterpart by the real structure.

From the construction in Section \ref{s:ps}, in order to eliminate
the base locus of $|(n-2)F|$, it is enough to 
eliminate the base locus of the system $|\mathscr L_1|$ on $Z_1$,
where $\mathscr L_1$ is the line bundle \eqref{L}.
We have the morphism $f_1:Z_1\to \Lambda\simeq\CP^1$ whose fibers are strict transforms
of the members of the pencil $|F|$.
Also recall that the restriction of $\mathscr L_1$ to
a general fiber $S$ of $f_1$ is isomorphic to the line bundle $(n-2)K_S^{-1}$
with the fixed components \eqref{bc0} removed.
Let 
\begin{align}
\lambda_1,\lambda_2,\cdots,\lambda_{n-1}
\end{align}
be points on the rational normal curve $\Lambda$ which correspond
to the reducible members $S_i^+ + S_i^-\in |F|$,
$1\le i\le n-1$, respectively.
Of course, the collection $f_1^{-1}(\lmd_i)$, $1\le i\le n-1$, are all reducible fibers of $f_1$.
Recall that as defined in \eqref{cij},
on the divisors $S_i^+$ and $S_i^-$ in $Z_1$,
there is a curve $C_{i,j}$ which is identified with
the curve $C_j$ in $Z$ under the birational morphism
$\mu:Z_1\to Z$.

The following property of the base locus of the system $|\mathscr L_1|$ 
follows immediately from Lemma \ref{lemma:int1}.
(See also the right picture in Figure \ref{fig:L1},
where the base curves lying on $S_{n-2}^+\cup S_{n-2}^-$
are written as segments with small triangles in the case 
$n=7$.)
\begin{proposition}\label{prop:bc8}
The  base locus of $|\mathscr L_1|$ contains the following curves:
\begin{align}\label{bc8}
\bigsqcup_{3\le i\le n-2} 
\bigg(
\bigcup_{3\le j\le i} C_{i,j}
\bigg)
\quad {\text{and}} \quad
\bigsqcup_{3\le i\le n-2} 
\bigg(
\bigcup_{3\le j\le i} \ol C_{i,j}
\bigg),
\end{align}
and 
\begin{align}\label{bc9}
C_{n-1,1}\sqcup \ol C_{n-1,1}.
\end{align}
\end{proposition}

Note that for each $i$ with $3\le i\le n-2$ the curves
in the parentheses of  \eqref{bc8}
are connected. So
each of the two curves \eqref{bc8} consists of $(n-4)$ connected components.
In particular they are empty when $n=4$.

\begin{remark}
{\em
As we shall see below, the two base curves
\eqref{bc8} and \eqref{bc9} have quite different nature.
}
\end{remark}
 
\begin{figure}
\includegraphics{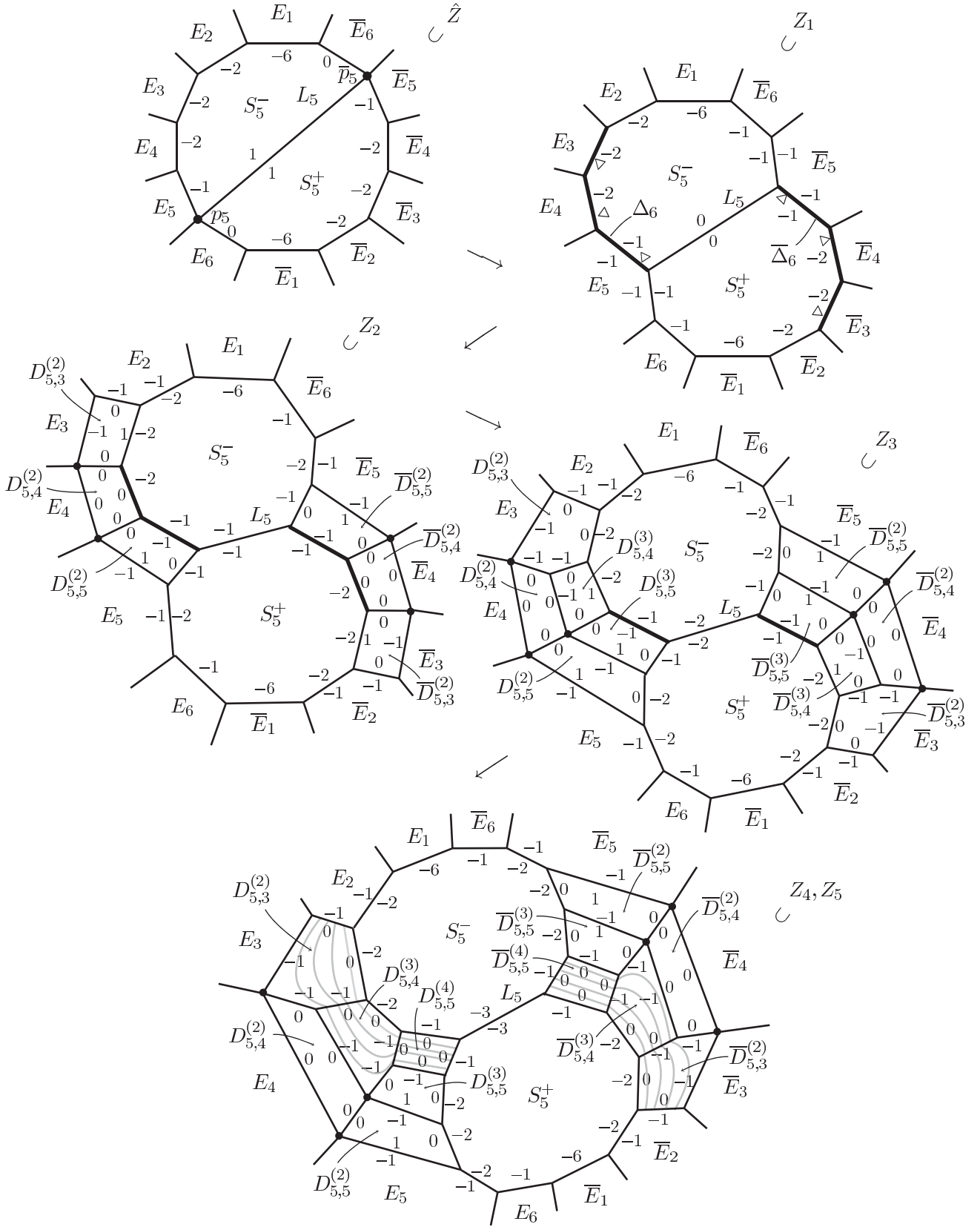}
\caption{
The sequence of blowups  in the case $n=7$, 
which eliminates the base locus of the system
$|\mathscr L_{n-2}|$ lying on $S_5^+\cup S_5^-$.}
\label{fig:elim1}
\end{figure}

\begin{figure}
\includegraphics{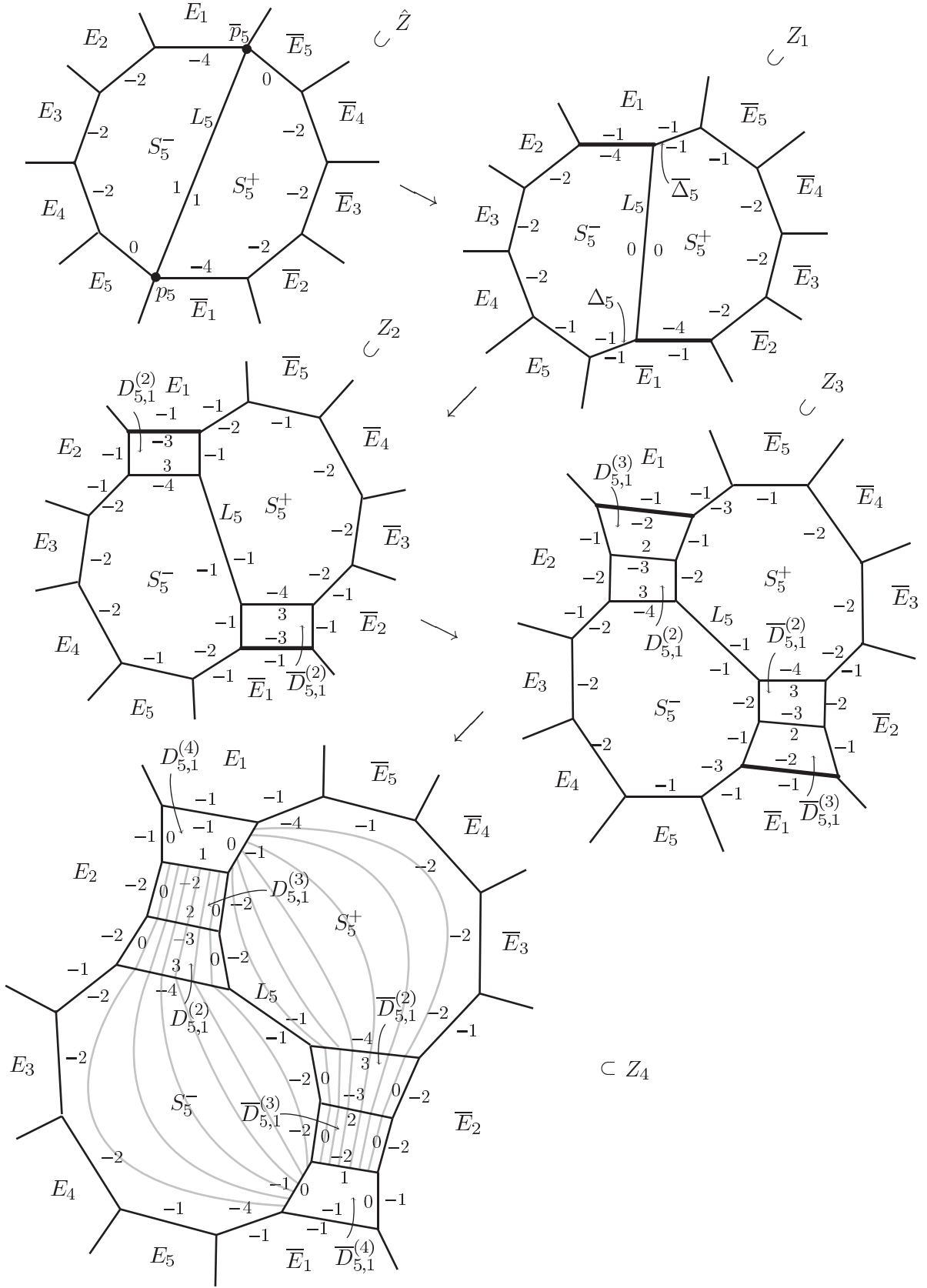}
\caption{
The sequence of blowups  in the case $n=6$, 
which eliminates the base locus of the system
$|\mathscr L_{n-2}|$ lying on $S_5^+\cup S_5^-$.}
\label{fig:elim2}
\end{figure}

Let $\mu_2:Z_2\to Z_1$ be the blowup of $Z_1$ at 
all the curves \eqref{bc8} and \eqref{bc9},
and $D^{(2)}_{i,j}$  the exceptional divisor over
$C_{i,j}$.
Each $D^{(2)}_{i,j}$ 
is isomorphic to a ruled surface
over $C_{i,j}$
(which is not isomorphic to $\CP^1\times\CP^1$ in general).
When $n>5$, as each connected component
of the center \eqref{bc8} is reducible
for $4\le i\le n-2$,
the variety $Z_2$ has an ODP at
the common point of the four divisors 
$$
E_{j},\,E_{j+1},\, D^{(2)}_{i,j},\, D^{(2)}_{i,j+1},
$$
where $4\le i\le n-2$ and $3\le j\le i-1$.
(In Figure \ref{fig:elim1}, 
these are indicated by dotted points 
in the case $(n,i) =(7,5)$.)
In particular, these ODP-s are not lying on  (the
strict transforms of) the divisors $S_i^+$ and $S_i^-$.
Define a line bundle $\mathscr L_2$ over $Z_2$ by
\begin{align}\label{L2}
\mathscr L_2:=\mu_2^*\mathscr L_1 - 
\sum_{i=3}^{n-2}
\sum_{j=3}^{i}
\big(D^{(2)}_{i,j} + \ol D^{(2)}_{i,j}\big)
-D^{(2)}_{n-1,1} - \ol D^{(2)}_{n-1,1}.
\end{align}
Namely we are just subtracting all the exceptional divisors from the pulled-back bundle.
The bundle $\mathscr L_2$ is clearly real.
Further, for each $i$ and $j$ with 
$3\le i\le n-2$ and $3\le j\le i$, we define
 curves on $Z_2$ by
\begin{align}\label{bc10}
C_{i,j}^{(2)}:= D^{(2)}_{i,j}\cap S_i^-
\quad{\text{and}}\quad
\ol C_{i,j}^{(2)}:= \ol D^{(2)}_{i,j}\cap S_i^+.
\end{align}
Since the variety $Z_2$ is non-singular 
on these curves, their intersection numbers
with the line bundle $\mathscr L_2$ make sense.
These can be computed in a similar way to 
the proof of Lemma \ref{lemma:int1}, and we obtain 
\begin{equation}
\label{int13}
\big(\mathscr L_2,  C_{i,j}^{(2)} \big)_{Z_2} = 
\left\{
\begin{array}{cc}
1, & j=3,\\
0, &  3<j<i,\\
-1, & j=i.
\end{array}
\right.
\end{equation}
Hence disposing the curves $C_{i,3}^{(2)}$
and $\ol C_{i,3}^{(2)}$ for each $i$, the curves
\begin{align}\label{bc30}
\bigsqcup_{4\le i\le n-2} 
\bigg(
\bigcup_{4\le j\le i} C_{i,j}^{(2)}
\bigg)
\quad {\text{and}} \quad
\bigsqcup_{4\le i\le n-2} 
\bigg(
\bigcup_{4\le j\le i} \ol C_{i,j}^{(2)}
\bigg)
\end{align}
are base curves of the system $|\mathscr L_2|$.
(In Figure \ref{fig:elim1}, these curves are indicated
by bold segments in the case $(n,i) =(7,5)$.)
In particular, not only the number of the base curves
lying on $S_i^-$ for fixed $i$,
but also the number of the connected components
of the base curves decrease by one as an effect of the 
blowup $\mu_2$, and
 \eqref{bc30} is empty when $n=5$.

Next, in order to inspect the base locus on the 
isolated exceptional divisor 
$D_{n-1,1}^{(2)}$, we put
\begin{align}\label{bc11}
C^{(2)}_{n-1,1}:= D_{n-1,1}^{(2)}\cap E_1
\quad{\text{and}}\quad
\ol C^{(2)}_{n-1,1}:=\ol D_{n-1,1}^{(2)}\cap \ol E_1.
\end{align}
(In Figure \ref{fig:elim2} these curves are indicated by
bold segments in the case $n=6$.)
Here note that unlike the above  curves \eqref{bc10},
we are taking intersection with the exceptional divisors $E_1$ and $\ol E_1$,
and hence
these curves are not lying on $S_i^-$ nor $S_i^+$.
Then $Z_2$ is non-singular on these curves,
and we can compute as
\begin{align}\label{incr}
\big(\mathscr L_2,  C_{n-1,1}^{(2)} \big)_{Z_2} = 
4-n.
\end{align}
In view of Lemma \ref{lemma:int1},
this increases from $(\mathscr L_1, C_{n-1,1})_{Z_1}$
by one. 
In particular, the curves \eqref{bc11} are base curves
when $n>4$.

When $n=4$, we finish the operation here
(i.e.\,at $Z_2$).
If $n>4$, let $\mu_3:Z_3\to Z_2$ be the blowup at
the curves \eqref{bc30} and \eqref{bc11}.
Let $D^{(3)}_{i,j}\subset Z_3$ be the exceptional divisor over the curve $C_{i,j}^{(2)}$,
and define a line bundle over $Z_3$ by
\begin{align}\label{L3}
\mathscr L_3:=\mu_2^*\mathscr L_2 - 
\sum_{i=4}^{n-2}
\sum_{j=4}^{i}
\big(D^{(3)}_{i,j} + \ol D^{(3)}_{i,j}\big)
-D^{(3)}_{n-1,1} - \ol D^{(3)}_{n-1,1},
\end{align}
which is still real.
Then by the same reason for the blowup $\mu_2:Z_2\to Z_1$, if $n>6$, the variety $Z_3$ has new ordinary double points
 over the singular points of the curves \eqref{bc30}.
But again they are not lying on (the strict transform of) the divisors $S_i^+$ and $S_i^-$.
If we put
\begin{align}\label{bc12}
C_{i,j}^{(3)}:= D^{(3)}_{i,j}\cap S_i^-
\quad{\text{and}}\quad
\ol C_{i,j}^{(3)}:= \ol D^{(3)}_{i,j}\cap S_i^+,
\end{align}
then again by computing intersection numbers with $\mathscr L_3$, we deduce that these are base curves of $|\mathscr L_3|$ as long as $5\le i\le n-2$ and $5\le j\le i$.
(In Figure \ref{fig:elim1} these curves are indicated by bold
segments in the case $(n,i) =(7,5)$.)
In particular, the number of the base curves
lying on $S_i^-$ for fixed $i$ and 
the number of the connected components
of the base curves
again decrease by one by the effect
of the blowup $\mu_3$.
Also, the intersection numbers of $\mathscr L_3$
with the curves
\begin{align}\label{bc13}
C^{(3)}_{n-1,1}:= D_{n-1,1}^{(3)}\cap E_1
\quad{\text{and}}\quad
C^{(3)}_{n-1,1}:=\ol D_{n-1,1}^{(3)}\cap \ol E_1
\end{align}
increase by one from \eqref{incr},
and the curves \eqref{bc13} are base curves if $n>5$.
(In Figure \ref{fig:elim2} these curves are 
indicated by bold segments in the case $n=6$.)
When $n=5$, we stop the operation here.
If $n>5$, we blowup $Z_3$ at the base curves \eqref{bc12} and \eqref{bc13},
and repeat the same operation above.

Continuing this process,
for the twistor space on $n\CP^2$
with $n$ being arbitrary, we obtain a sequence of the explicit blowups
$$
Z_{n-2} \stackrel{\mu_{n-2}}{\lra} Z_{n-3} 
\stackrel{\mu_{n-3}}{\lra} \cdots
\stackrel{\mu_{3}}{\lra} Z_2 
\stackrel{\mu_{2}}{\lra} Z_1.
$$
together with a real line bundle $\mathscr L_m$ over each $Z_m$.
Let $\nu:Z_{n-2}\to Z_1$ be the composition of all these blowups.
(We note that the number of times of blowups are $(n-4)$
for the base curves \eqref{bc8}, and 
$(n-3)$ for the base curves \eqref{bc9}.
The final blowup $\mu_{n-2}:Z_{n-2}\to Z_{n-3}$
changes only the latter base curves \eqref{bc9}
and does not change the former base curve \eqref{bc8}.)
Then  the exceptional divisor of $\nu$ consists of
\begin{align}\label{ed1}
\nu^{-1} (C_{i,j}) = \bigcup_{k=2}^{j-1 } D^{(k)}_{i,j},\quad
{\text{where}}\quad 3\le i\le n-2
\quad
{\text{and}}\quad 3\le j\le i,
\end{align}
\begin{align}\label{ed2}
\nu^{-1} (C_{n-1, 1}) = \bigcup_{k=2}^{ n-2 } D^{(k)}_{n-1,1}
\end{align}
and the images of these divisors by the real structure.
From the explicit construction, for divisors in \eqref{ed2} we readily have
\begin{align}\label{ed3}
D^{(k)}_{n-1,1}\simeq \Sigma_{n-k-1}
\end{align}
where $\Sigma_d$ denotes the ruled surface of degree $d$ over $\CP^1$.
In particular, the divisor $D_{n-1,1}^{(n-2)}$ (obtained 
from the final blowup) is isomorphic to 
one point blown-up of $\CP^2$.
Further, the intersection of the adjacent components
$D^{(k)}_{n-1,1}$ and 
$D^{(k+1)}_{n-1,1}$ is always a section of the ruling.
Furthermore, these sections are mapped isomorphically
to the curve $C_{n-1,1}$ by $\nu$.
Thus it would be possible to say that
 the divisor \eqref{ed2} has a structure of 
a ladder over $C_{n-1,1}$ by $\nu$.
Similarly, for each $i$ and $j$, the divisor \eqref{ed1} forms
a ladder over $C_{i,j}$,
but this ladder is growing up in the reverse direction
with the ladder \eqref{ed2}.

Then we have the following

\begin{proposition}\label{prop:elim}
The linear system $|\mathscr L_{n-2}|$ on the variety $Z_{n-2}$ obtained above
is base point free.
\end{proposition}

\proof
Let $Z_{n-2}' \to Z_{n-2}$ be any small deformation of all ODP-s on $Z_{n-2}$
which preserves the real structure,
and let $\mathscr L'_{n-2}$ be the pullback of $\mathscr L_{n-2}$ to $Z_{n-2}$.
It suffices to show that $|\mathscr L'_{n-2}|$ is base point free.
For this purpose, we compute
the restrictions of the line bundle $\mathscr L'_{n-2}$
to the exceptional divisors of $\nu$ and $\mu$
(namely $D_{i,j}^{(k)}$ and $E_j$).
In the sequel we obtain surjective morphisms
from each of these divisors to $\CP^1$ (including
the above ruling map for $D_{n-1,1}^{(k)}\simeq\Sigma_{n-k-1}$),
for which we use the common letter $p$. 
(In Figures \ref{fig:elim1} and \ref{fig:elim2},   fibers of these morphisms are
indicated by gray curves in the cases 
$(n,i) =(7,5)$ and $(n,i)=(6,5)$ respectively.)

First for the exceptional divisor $D_{i,j}^{(k)}$ in \eqref{ed1},
 by computing the intersection numbers of $\mathscr L'_{n-2}$ with
 curves in $D_{i,j}^{(k)}$ which are obtained as an intersection 
of other exceptional divisors,
it is possible to show that the line bundle 
 $\mathscr L'_{n-2}$ is {\em trivial} over the  divisors
\begin{align}\label{ed0}
D_{i,j}^{(k)}, \,\,3\le i\le n-2, \,\,3\le j\le i,
\,\,2\le k\le j-2.
\end{align} 
Among all the divisors \eqref{ed1}, these are characterized by  disjointness
with the divisor $S_i^-$
(i.e.\,we are just excluding the case $k=j-1$ from
the divisors \eqref{ed1}.) 
Note that for any fixed $i$,
the union of all the divisors \eqref{ed0} is connected.
On the other hand, over the remaining exceptional divisors
\begin{align}\label{ed4}
D_{i,j}^{(j-1)},\,\,3\le i\le n-2,\,\,3\le j\le i,
\end{align}
the line bundle
$\mathscr L'_{n-2}$ is not trivial but  of the form $p^*\mathscr O(1)$,
where $p$ is a surjective morphism $D_{i,j}^{(j-1)}\to\CP^1$ which has the intersection curve
$D_{i,j}^{(j-1)}\cap S_i^-$ as a smooth fiber.

Next for the restriction to the divisor $D^{(k)}_{n-1,1}$ in \eqref{ed2} (or \eqref{ed3}), for any $k$ with $2\le k<n-2$, 
over $D^{(k)}_{n-1,1}$, the line bundle $\mathscr L_{n-2}'$ is 
isomorphic to $p^*\mathscr O(1)$, where $p:
D^{(k)}_{n-1,1}\to\CP^1$ is the ruling map.
Over the remaining divisor $D_{n-1,1}^{(n-2)}$
(which is placed at an end of the ladder), 
$\mathscr L_{n-2}'$ is of the form $\varepsilon ^*\mathscr O(1)$,
where $\varepsilon: D_{n-1,1}^{(n-2)}\simeq\Sigma_1\to \CP^2$ is a blow-down.

For the restriction to the exceptional divisor $E_i$ ($1\le i\le n-1$),
we recall that the line bundle $\mathscr L_1$ over 
$Z_1$ is trivial over $E_{n-1}$ by
Remark \ref{rmk:triv}.
On the other hand, $\mathscr L_1$ is non-trivial over
the remanning divisors $E_i$, $i<n-1$.
However, as an effect of the blowups in $\nu$,
the final line bundle
$\mathscr L_{n-2}'$ is {\em trivial} over any divisors $E_i$,
so far as $i\neq 2$.
On the other hand, over the component $E_2$, $\mathscr L_{n-2}'$ is of the form
$p^*\mathscr O(1)$, where $p:E_2\to \CP^1$ is a surjective morphism
for which  the intersection $S\cap E_2$ is a section,
with $S$ being any fiber of the composition
\begin{align}\label{f00}
Z'_{n-2} \lra Z_{n-2} \stackrel{\nu}{\lra}
Z_1 \stackrel{f_1}{\lra} \CP^1.
\end{align}

We also need to know the restriction of $\mathscr L_{n-2}'$ to
the divisor $S_{n-1}^-$.
For this we again notice that there is a surjective morphism $p:S_{n-1}^-\to\CP^1$
for which the two intersection curves $S_{n-1}^-\cap D_{n-1,1}^{(2)}$ 
and $S_{n-1}^-\cap \ol D_{n-1,1}^{(n-2)}$ are (mutually disjoint) sections of the morphism.
(The morphism is induced by the linear system
$|\sum_{i=2}^{n-1} C_{n-1,i} + \ol C_{n-1,1}|$
on $S_{n-1}^-$.)
Then the restriction of $\mathscr L'_{n-2}$
to $S_{n-1}^-$ is again of the form $p^*\mathscr O(1)$.

Utilizing all these restriction data, we show that $|\mathscr L_{n-2}'|$ is base point free.
First it is clear from the beginning that 
$\Bs\,|\mathscr L_{n-2}'|$ is contained in the exceptional divisors
 \eqref{ed1}, \eqref{ed2}, $E_i$ $(1\le i\le n-1)$, or their conjugations
by the real structure.
If there is a base point on the divisors \eqref{ed0}
for some $i$ (so that $3\le j\le i\le n-2$), then by the triviality of $\mathscr L'_{n-2}$
over these divisors,
all the divisors \eqref{ed0} must be  fixed components
of $|\mathscr L'_{n-2}|$ for the above $i$.
Then since the component $D_{i,i}^{(2)}$ intersects $E_{i}$ and
$\mathscr L'_{n-2}$ is trivial over 
$\cup_{j\neq2} (E_j\cup\ol E_j)$ as above,
whole of the last union must also be fixed components of 
$|\mathscr L_{n-2}'|$.
This means that the original system $|\mathscr L_1|$ on the space $Z_1$ has
$E_i$ $(i\neq 2)$ as a fixed component.
But this contradicts the fact that
$|\mathscr L_1|$ does not have a base point
on smooth fibers of $f_1:Z_1\to\CP^1$.
Therefore $\Bs\,|\mathscr L_{n-2}'|$ is disjoint from 
the divisors \eqref{ed0} and also from their conjugations.
Moreover if there is a base point on the divisor $D_{i,j}^{(j-1)}$ in \eqref{ed4}, 
then through the  morphism $p:D_{i,j}^{(j-1)}\to \CP^1$ (where we decrease $j$ one by one until it becomes $3$),
we finally obtain that there must be a base point on $E_2$.
Hence $\Bs\,|\mathscr L_{n-2}'|\cap E_2$ has to be a fiber of the
above morphism $p:E_2\to \CP^1$.
But since fibers of the last morphism intersect
any fiber of the morphism \eqref{f00},
this again contradicts the fact that $|\mathscr L_1|$ does not have a base point on 
smooth fibers of $f_1$.
So we conclude $\Bs\,|\mathscr L_{n-2}'|\cap D_{i,j}^{(j-1)}=\emptyset$ also for the divisors in \eqref{ed4}.
Furthermore, these argument clearly imply that
$\Bs\,|\mathscr L_{n-2}'|\cap E_i = \emptyset$
for any $1\le i\le n-1$.

It remains to see that $|\mathscr L_{n-2}'|$ does not have
a base point on the ladder \eqref{ed2}.
If the system $|\mathscr L_{n-2}'|$ has a base point on the divisor 
$D^{(k)}_{n-1,1}$ for some $k<n-2$, then 
via the ruling map $D^{(k)}_{n-1,1}\to\CP^1$
(where we decrease $k$ until it becomes $2$),
the system has a base point on the intersection curve $D^{(2)}_{n-1,1}\cap S_{n-1}^-$,
which is a section of the morphism
$p:S_{n-1}^-\to \CP^1$.
Hence $|\mathscr L_{n-2}'|$ has a base point along a fiber of the  morphism
$p:S_{n-1}^-\to \CP^1$.
But this contradicts the fact that $\Bs\,|\mathscr L_{n-2}'|$ is contained in the 
divisor \eqref{ed2} (which is already proved).
Hence if $k<n-2$ we have $ D^{(k)}_{n-1,1} \cap \Bs\,|\mathscr L_{n-2}'| =\emptyset$.
Finally we show $ D^{(n-2)}_{n-1,1} \cap \Bs\,|\mathscr L_{n-2}'| =\emptyset$.
Let $\Phi'_{n-2}:Z'_{n-2}\to\CP^n$ be the rational map
associated to $|\mathscr L_{n-2}'|$.
Then since the restriction of
$\mathscr L'_{n-2}$ to the divisor
$D_{n-1,1}^{(k)}$ is of the form $p^*\mathscr O(1)$
and also that 
$\Bs\,|\mathscr L_{n-2}'| \cap D_{n-1,1}^{(k)}
=\emptyset$ for $2\le k < n-2$,
the image $\Phi'_{n-2} (D_{n-1,1}^{(k)})$ must be
a line for this range of $k$.
Moreover this line is independent of $k$,
because of the ladder structure of 
the divisor \eqref{ed2}.
The ladder structure also implies that
the curve $S_{n-1}^-\cap D_{n-1,1}^{(2)}$
(which is a section of the morphism $p:S_{n-1}^-\to\CP^1$)
and the curve 
$D_{n-1,1}^{(n-2)} \cap D_{n-1,1}^{(n-3)}$ 
are also mapped to the same line by $\Phi'_{n-2}$.
Hence, via the morphism $p:S_{n-1}^-\to\CP^1$,
the other section $S_{n-1}^-\cap \ol D_{n-1,1}^{(n-2)}$ is 
also mapped to the same line by $\Phi'_{n-2}$.
Therefore, by the real structure, the image 
 $\Phi'_{n-2} ( S_{n-1}^+\cap D_{n-1,1}^{(n-2)})$ of the conjugate curve is a line.
 Hence $\Phi'_{n-2}$ maps
 the two curves 
$ D^{(n-2)}_{n-1,1} \cap  D^{(n-3)}_{n-1,1}$  and $ S_{n-1}^+\cap  D_{n-1,1}^{(n-2)}$ on the surface $D^{(n-2)}_{n-1,1}\simeq\Sigma_1$ 
to lines.
However, these lines cannot be identical, since 
the former curve belongs to $|\varepsilon^*\mathscr O(1)|$
on $D^{(n-2)}_{n-1,1}$, while the latter does not.
Hence we conclude that the image $\Phi'_{n-2}(D^{(n-2)}_{n-1,1}) $ contains two lines.
Therefore $\Phi'_{n-2}(D^{(n-2)}_{n-1,1}) $ has to
be two-dimensional, meaning that 
$\Phi'_{n-2}(D^{(n-2)}_{n-1,1})  = \pi^{-1}(\lmd_{n-1}).$
Hence recalling that $\mathscr L'_{n-2}$
is isomorphic to 
$\varepsilon^*\mathscr O_{\CP^2}(1)$
over $ D^{(n-2)}_{n-1,1}$, 
we obtain $ D^{(n-2)}_{n-1,1} \cap \Bs\,|\mathscr L_{n-2}'| =\emptyset$.

Thus we have completed a proof of 
Proposition \ref{prop:elim}.
\proofend

\vsp
By using this elimination, we prove the following 
result concerning the behavior of the rational
map $\Phi$ on the degree-one divisors,
which will be needed in the next section.
Recall that for $1\le i\le n-1$, $L_i$ means the twistor line $S_i^+\cap S_i^-$.
\begin{proposition}\label{prop:image4}
Let $\Phi:Z\to Y\subset\CP^n$ be the rational map associated
to the system $|(n-2)F|$ as before.
\begin{enumerate}
\item[(i)] If $1\le i<n-1$, both of the restrictions 
$\Phi|_{S_i^+}$ and $\Phi|_{S_i^-}$ are birational over
the plane $\pi^{-1}(\lmd_i)$.
Moreover, the image $\Phi(L_i)$ is a conic in the plane. 
\item[(ii)] Both of the images $\Phi(S_{n-1}^+)$ and $\Phi(S_{n-1}^-)$ are lines in the plane $\pi^{-1}(\lmd_{n-1})$.
Further, these lines are distinct.
\end{enumerate}
\end{proposition}

\proof
From our construction, it is enough to show the same claim for
the morphism $\Phi'_{n-2}$ associated to $|\mathscr L'_{n-2}|$ on $Z'_{n-2}$.
For the birationality in (i), as we know that $\Phi'_{n-2}$ is of degree two preserving the
real structure, it is enough
to show that $\Phi'_{n-2}|_{S_i^-}$ is surjective over the plane
$\pi^{-1}(\lmd_i)$.
As in the  above proof of
Proposition \ref{prop:elim}, 
the line bundle $\mathscr L'_{n-2}$ over $Z'_{n-2}$ is
trivial over the exceptional divisors \eqref{ed0}, and therefore
their images by $\Phi'_{n-2}$ must be a point.
Similarly, since the restriction of $\mathscr L'_{n-2}$ 
to the divisor $D_{i,j}^{(j-1)}$ in \eqref{ed4}
is of the form $p^*\mathscr O(1)$ where $p:D_{i,j}^{(j-1)}\to
\CP^1$ is the surjective morphism,
the image  $\Phi'_{n-2}(D_{i,j}^{(j-1)})$ must be a line.
Therefore, since we already know that $|\mathscr L_{n-2}'|$ 
is base point free, the remaining divisor
$S_i^-$ has to be mapped surjectively to the plane $\pi^{-1}
(\lmd_{n-1})$, as claimed.

For the assertion about the images of the twistor lines in (i), since
$\Bs\,|\mathscr L_{n-2}'|=\emptyset$ and
$\Phi'_{n-2}$ is degree one over $S_i^-$
as above,
it suffices to show that 
$(\mathscr L_{n-2}', L_i)_{Z'_{n-2}}=2$
 for $1\le i<n-1$.
From the formula \eqref{L1-1}, 
on the manifold $Z_1$ we have
$$
(\mathscr L_1, L_i)_{Z_1} = 2(i-1),
$$
since $(f_1^*\mathscr O(1), L_i)_{Z_1} = 0$
(as $L_i\subset f_1^{-1}(\lmd_i)$) 
and, among the divisors $E_j$ and $\ol E_j$, $1\le j\le n-1$, only $E_i$ and $\ol E_i$ intersect $L_i$
and the intersections are transversal.
Then for each of the blowup $\mu_m:Z_m\to Z_{m-1}$, 
since we are removing the exceptional divisors of $\mu_m$ by
multiplicity one (see \eqref{L2} and \eqref{L3}),
if the center of the blowup intersects $L_i$,
the intersection number satisfies 
\begin{align}\label{dec1}
(\mathscr L_m,L_i)_{Z_m} = (\mathscr L_{m-1},L_i)_{Z_{m-1}}-2.
\end{align}
Further from the explicit centers of the blowups $\mu_m$, this actually happens exactly when $2\le m\le i-1$. 
Hence the decreasing \eqref{dec1} happens precisely
$(i-2)$ times.
Therefore we have $$(\mathscr L'_{n-2},L_i)_{Z_{n-2}'}
=(\mathscr L_{n-2},L_i)_{Z_{n-2}}
= (\mathscr L_1,L_i)_{Z_1} - 2(i-2) = 2,
$$ 
as claimed.

The second assertion (ii) is already shown in the 
final part of the above proof of
Proposition \ref{prop:elim}.
\proofend 

\begin{remark}
{\em
As showed in the final part of the proof of
Proposition \ref{prop:elim}, the final exceptional divisors $D^{(n-2)}_{n-1,1}$
and $\ol D^{(n-2)}_{n-1,1}$ are mapped 
birationally onto the plane $\pi^{-1}(\lmd_{n-1})$. 
}
\end{remark}

\section{Defining equation of the quartic hypersurface}
In this section, assembling all the results obtained so far,
we shall obtain defining equation of a
quartic hypersurface in $\CP^n$ which cut out
the branch divisor 
of the pluri-half-anticanonical map $\Phi:Z\to Y\subset\CP^n$.

\subsection{Double curves on the branch divisor, and a quadratic hypersurface
containing them}
\label{ss:dc}
In the study of plane quartic curves, the notion of bitangent has been played
a significant role.
Similarly, for quartic surfaces  in 
a projective space, a plane which touches the surface along a curve
is meaningful, because 
such a plane brings much information about a defining equation
of the surface.
This is also the case for the study of the 
branch divisor for the present twistor spaces.

As before let $Z$ be the twistor space on $n\CP^2$ which has 
the surface $S$ constructed in Section \ref{ss:S} as a real member of 
the system $|F|$, and $\Phi:Z\to \CP^n$ be the rational map
associated to the system $|(n-2)F|$.
If $\pi:\CP^n\to \CP^{n-2}$ denotes the natural linear projection 
corresponding to the subspace $S^{n-2} H^0(F)\subset H^0((n-2)F)$ as before,
the image $\Phi(Z)$ is the scroll $Y=\pi^{-1}(\Lambda)$,
where $\Lambda$ is a rational normal curve in $\CP^{n-2}$ (Proposition \ref{prop:Z2}).
We also know that the branch divisor $B$ of $\Phi:Z\to Y$ is 
of the form $Y\cap \mathscr B$, where $\mathscr B$ is 
a quartic hypersurface in $\CP^n$.
In order to determine a defining equation of $\mathscr B$,
we call  a curve on the branch divisor $B$ to be a {\em double curve}
if there exists a hyperplane $H\subset\CP^n$ such that 
$H\cap B$ is a non-reduced curve on the surface $B$.
Geometrically, this means that the hyperplane section $H\cap Y$ 
is tangent to $B$ along the curve.

As before let 
$
\lambda_1,\lambda_2,\cdots,\lambda_{n-1}
$
be points on the curve $\Lambda\subset \CP^{n-2}$ which correspond
to the reducible divisors $S_i^+ + S_i^-\in |F|$ respectively,
where we are identifying the curve $\Lambda$ with the 
space $\mathbb P^*H^0(F)\simeq\CP^1$ through the diagram
\eqref{016}.
By Proposition \ref{prop:image4} (i),
for $1\le i\le n-2$,
we have 
$\Phi(S_i^+) = \Phi(S_i^-) = \pi^{-1}(\lambda_i)$
and the image $\Phi(L_i)$ is a conic
in the plane $\pi^{-1}(\lmd_i)$.
We put 
\begin{align}
\mathscr C_i : = \Phi(L_i),
\quad 
1\le i\le n-2 
\end{align} 
for these conics.
Since $\Phi$ preserves the real structure,
these conics  are  real.
Then since $\Phi^{-1}(\pi^{-1}(\lambda_i))$ splits into the union $S_i^+\cup S_i^-$,
the plane $\pi^{-1}(\lambda_i)$ must touch the divisor $B$ 
along the conic $\mathscr C_i$.
Hence for $1\le i\le n-2$, the conic $\mathscr C_i$ is a double curve in the above sense.
We call these curves as {\em double conics}.
For the case $i=n-1$, by Proposition \ref{prop:image4} (ii),
the image $\Phi(S_i^+)$ and $\Phi(S_i^-)$ are mutually distinct lines.
We denote the union of these two lines
by  $\mathscr C_{n-1}$. Namely, we put
$$
\mathscr C_{n-1}:= \Phi(S_{n-1}^+\cup S_{n-1}^-).
$$
As in the case $i<n-1$, this also has to be a double curve on $B$,
and we call it as a  {\em splitting double conic}.

By Proposition \ref{prop:red}, there exist
two special reducible members of the system 
$|(n-2)F|$, both of which consist of
two irreducible components.
Let $H_n$ and $H_{n+1}$ be the hyperplanes in $\CP^n$ 
which correspond to these two reducible members.
Then since the preimage $\Phi^{-1}(H_n)$ (resp.\,$\Phi^{-1}(H_{n+1})$) 
 splits into the two irreducible components,
 if we denote $\mathscr C_n$ (resp.\,$\mathscr C_{n+1}$) for
the image of the intersection of the two irreducible components of 
$\Phi^{-1}(H_n)$ (resp.\,$\Phi^{-1}(H_{n+1})$),
these are also double curves on $B$.
If we write $l$ for the line which is the 
indeterminacy locus (the center) of the 
linear projection $\pi$,
the hyperplane section $H_i\cap Y$ $(i\in\{n, n+1\})$ is a 
cone over the curve $\Lambda$ whose vertex is the point $l\cap H_i$.
Moreover since the twice of the double curve
$\mathscr C_i$ ($i\in\{n,n+1\}$) is an element
of the system $|\mathscr O(4)|$ on the cone,
the curve $\mathscr C_i$ belongs to the  system $|\mathscr O(2)|$ on the cone.
In particular the degree of $\mathscr C_n$
and $\mathscr C_{n+1}$ in $\CP^n$ is $2(n-2)$.

Thus we have $(n-1)$ double conics 
$\mathscr C_1,\cdots, \mathscr C_{n-1}$
on the planes and
two double curves 
$\mathscr C_n,\mathscr C_{n+1}$ on the cones.
These double curves play an essential role for obtaining a defining equation of the
quartic hypersurface.

\begin{proposition}\label{prop:Q}
Let  $\mathscr C_1,\cdots,\mathscr C_{n+1}$ be the double curves on 
the branch divisor $B\subset Y$ as above.
Then there exists a quadratic hypersurface $Q$ in $\CP^n$ which contains 
all these double curves, and which is different from the scroll $Y$.
\end{proposition}

\proof
For $1\le i\le n-1$ we denote $P_i$ for the plane $\pi^{-1}(\lmd_i)$.
We also write $D_n := Y\cap H_n$ and $D_{n+1} : = Y\cap H_{n+1}$
for the cones, on which $\mathscr C_n$ and $\mathscr C_{n+1}$ lie
respectively.
Let $\delta:\tilde Y\to Y$ be the blowup at the line (ridge) $l$,
and $\Sigma$ the exceptional divisor.
$\Sigma$ is biholomorphic to $\CP^1\times \CP^1$.
Let $\tilde{\pi}:\tilde Y\to \CP^1$ be the composition $\pi\circ\delta$,
which is clearly a morphism.

We again use the same letters for the strict transforms into $\tilde Y$ of the divisors and curves in $Y$.
On the resolved space $\tilde Y$, 
the divisors $D_n$ and $D_{n+1}$ are smooth and biholomorphic to $\Sigma_{n-2}$,
the ruled surface of degree $(n-2)$.
Moreover, on $\tilde Y$, the union
\begin{align}\label{rest20}
D:=\bigg(
\bigcup_{i=1}^{ n-1} P_i
\bigg)
\cup
D_n \cup D_{n+1} 
\cup \Sigma
\end{align}
is not only smooth normal crossing but also simply connected.
(Note that the fundamental group of 
 the union $D_n\cup D_{n+1}\cup \Sigma$ is  $\mathbb Z$,
but the generator becomes homotopic to the identity after adding the plane $P_i$.)

As the degree of $\Lambda$ in $\CP^{n-2}$ 
is $(n-2)$, we have
\begin{align}\label{rel1}
\delta^*\mathscr O_Y(1) \sim  \Sigma +  \tilde{\pi}^*\mathscr O_{\Lmd}(n-2). 
\end{align}
Therefore,
noting $P_i\in |\tilde{\pi}^*\mathscr O_{\Lmd}(1)|$ and $D_j\in |\delta^*\mathscr O_Y(1)|$
on $\tilde Y$,
we obtain that the divisor \eqref{rest20} (with all components counted just once) belongs to the 
linear system
$$
\big|
3 \Sigma +  \tilde{\pi}^*\mathscr O_{\Lmd}(3n-5) 
\big|.
$$
Hence again by using the relation \eqref{rel1}
 we obtain an exact sequence
\begin{align}
 \label{ses5}
0 \,\lra\,
- \Sigma +  \tilde{\pi}^*\mathscr O_{\Lmd}(1-n) 
\,\lra\, 
\delta^*\mathscr O_Y(2)
\,\lra\, 
\delta^*\mathscr O_Y(2)\big|_D
\,\lra\, 0.
\end{align}
Now a defining quartic polynomial of the hypersurface $\mathscr B$
 gives a section of the line bundle $\mathscr O_Y(4)$ as a restriction.
Let $s\in H^0(\mathscr O_Y(4))$ be this section. 
Then because all $\mathscr C_i$-s are double curves of $B$,
 for each $1\le i\le n-1$ and $j\in \{n,n+1\}$ we can take sections
 $$
 t_i\in H^0(P_i,
 \mathscr O_{P_i}(2)) \quad {\text{and}} \quad t_j\in H^0(D_j,\mathscr O_{D_j}(2))
$$
such that $  s|_{P_i} = t_i^2$ and
$  s|_{D_j} = t_j^2$.
Also, since the intersection $B\cap l$ consists of two points
(see Proposition \ref{prop:sing2}), 
there also exists an element $t_0
\in H^0(l,\mathscr O(2))$ such that 
$t_0^2= s|_l$.
 All these sections $t_i$, $t_j$ and $t_0$ are determined only up to sign.
Let $\tilde t_i,\,\tilde t_j$ and 
$\tilde t_0$ be the natural lifts 
of these sections by the blowup $\delta$.
Then from the simply connectedness of the divisor $D$, 
 we can choose the signs in a way that 
 any two agree on the intersections of 
the components of $D$.
Hence since
$D$ in $\tilde Y$ is smooth normal crossing, we obtain an element $t_D\in H^0(D,\delta^*\mathscr O_Y(2))$
such that it defines the double curves on
the plane $P_i$ and the cone $D_j$.

 On the other hand, by restricting the line bundle $\tilde{\pi}^*\mathscr O_{\Lmd}(1)$
 to the divisor $D_n$, it is easy to see that 
 $H^1 (\tilde Y, - \Sigma +  \tilde{\pi}^*\mathscr O_{\Lmd}(1-n) ) = 0$.
Hence from the exact sequence \eqref{ses5},
the restriction map 
$H^0(\delta^*\mathscr O_Y(2))
\to 
H^0(\delta^*\mathscr O_Y(2)\big|_D)$  is surjective.
Hence we obtain an element $\tilde t\in H^0(\delta^*\mathscr O_Y(2))$ such that 
$\tilde t|_D= t_D$.
Letting $t\in H^0(Y,\mathscr O_Y(2))$ be the element corresponding to $\tilde t$,
the divisor $(t)$ on $Y$ contains all the $(n+1)$ double curves.
Then since the restriction 
$H^0(\CP^n, \mathscr O(2))
\to H^0(Y, \mathscr O(2))$ is surjective,
 there is an element
$q\in H^0(\CP^n, \mathscr O(2))$ such that 
$q|_Y = t$.
Putting $Q:=(q)$, we obtain the required
quadratic hypersurface $Q$.
\proofend

\subsection{Defining equation of the 
quartic hypersurface}
\label{s:de}
Now we are ready for determining defining equation of the branch divisor.
Let all the notations be as in the last
subsection, and
we choose  homogeneous coordinates
$(z_0,z_1,\cdots,z_{n-2})$ on $\CP^{n-2}$ such that the rational normal curve $\Lmd$ is
realized as an image of the standard holomorphic map
\begin{align}\label{rnc}
\CP^1\ni (u_0,u_1)\mapsto 
(u_0^{n-2}, u_0^{n-3}u_1,\cdots, u_1^{n-2})
\in \CP^{n-2}.
\end{align}
We further normalize the coordinates 
in such a way that the last point $\lmd_{n-1}$  is the point $(0,0,\cdots,0,1)\in\CP^{n-2}$.
Then the hyperplane $\{z_0=0\}\subset\CP^{n-2}$ intersects $\Lmd$ at $\lmd_{n-1}$ by 
the highest multiplicity $(n-2)$. 
Under these normalizations of the coordinates,
our main result is stated as follows:
\begin{theorem}\label{thm:main}
Let $\Phi:Z\to Y\subset\CP^n$,  $\mathscr B$, $\Lmd$ and $(z_0,z_1,\cdots,z_{n-2})$ be as above.
Then for appropriate homogeneous coordinates $(z_0,z_1,\cdots,z_n)$ on $\CP^n$ which are
obtained as an extension of the above ones
on $\CP^{n-2}$,
the quartic hypersurface $\mathscr B$ is defined by the following polynomial:
\begin{align}\label{quart1}
z_0 z_{n-1}z_n f(z_0,z_1,\cdots,z_{n-2})= Q(z_0,z_1,\cdots,z_n)^2,
\end{align}
where $f(z_0,z_1,\cdots,z_{n-2})$ is a linear form  (not on $\CP^n$ but) on $\CP^{n-2}$, and
$Q(z_0,z_1,\cdots,z_n)$ is a quadratic form on $\CP^n$.
\end{theorem}

\begin{remark}{\em
As the proof below shows, the quadratic form $Q(z_0,\cdots,z_n)$ is exactly the defining equation
of the hyperquadric  in Proposition \ref{prop:Q}.
Therefore, since the restriction of the hyperquadric $Q$ to
the plane $\pi^{-1}(\lmd_{n-1}) = \{z_0=z_1=\cdots
=z_{n-3} = 0\}$ is 
the  the splitting double conic
$\mathscr C_{n-1}$,
we have a constraint that {\em the conic
defined by the equation
\begin{align}\label{conic5}
Q(0,0,\cdots,0,z_{n-2},z_{n-1},z_n) = 0
\end{align} 
is reducible} under the above normalization of
the coordinates.
}
\end{remark}

\begin{remark}
{\em
As the proof below shows, up to a non-zero constant, the linear form $f$ in the equation
\eqref{quart1} is uniquely determined from the $(n-2)$ points
$\lmd_1,\cdots,\lmd_{n-2}\in\Lmd$.
}
\end{remark}

\begin{remark}{\em
At first sight one might think that when $n=4$ the equation \eqref{quart1} coincides with 
the equation (1.2) in \cite{HonDS4_1} of the quartic hypersurface
in $\CP^4$.
But this is not correct, since the linear polynomial $f$ in \cite[(1.2)]{HonDS4_1} includes
not only $z_0,z_1,\cdots,z_{n-2}$ but also $z_{n-1}$ and $z_n$.
(This is not a minor difference, as the type of the singularities of the branch divisor 
becomes quite different.)
This means that the twistor spaces in this paper is not a direct generalization of the
twistor spaces studied in \cite{HonDS4_1}.
Rather, they are a direct generalization of one type of the twistor spaces  studied in \cite{HonDS4_2},
which we  call `type II' there.
}
\end{remark}

\noindent
{\em Proof of Theorem \ref{thm:main}.}
First  let $z_{n-1}$ and $z_n$ be  linear forms on $\CP^n$ 
such that $(z_{n-1}) = H_{n}$ and $(z_n) = H_{n+1}$,
where as before  $H_{n}$ and $H_{n+1}$ are the hyperplanes corresponding to
the two reducible members of $|(n-2)F|$.
Then $(z_0,z_1,\cdots,z_n)$ provides homogeneous coordinates on $\CP^n$,
with respect to which the line $l$ is defined by $z_0=z_1=\cdots=z_{n-2}=0$.

For an algebraic variety $X\subset\CP^n$,
we denote by $I_X\subset\mathbb C[z_0,\cdots,z_n]$
for the homogeneous ideal of $X$.
Let $F=F(z_0,\cdots,z_n)$ be a defining quartic
polynomial of $\mathscr B$.
Obviously $F$ is defined only up to an ideal
$I_{Y}\subset\mathbb C[z_0,\cdots,z_n]$.
Let $Q(z_0,\cdots,z_n)$ be a defining polynomial of the hyperquadric $Q$
whose existence was proved in Proposition \ref{prop:Q}.
$Q$ contains all the double curves $\mathscr C_i$, $1\le i\le n+1$.
For each $i$ with $1\le i\le n-1$ let $ P_i\subset Y$
be the plane $\pi^{-1}(\lambda_i)$ as before.
Then as $(F|_{ P_i}) = 2\mathscr C_i
= (Q^2|_{ P_i})$
as divisors on the plane $ P_i$,
there exists a constant $c_i$ such that $F-c_iQ^2 \in
I_{ P_i}\subset\CC[z_0,\cdots,z_n]$.
If $c_i\neq c_j$ for some $i\neq j$,
we obtain $Q^2\in I_{ P_i} + I_{ P_j}$.
Further the last ideal is readily seen to be equal to
$ I_{ P_i\cap P_j}$, and therefore
equals to $I_l=(z_0,z_1,\cdots,z_{n-2})\subset\CC[z_0,\cdots,z_n]$.
Hence $Q\in (z_0,z_1,\cdots,z_{n-2})$.
But this means that the divisor  $(Q|_{ P_i})$
contains the line $l$, which contradicts
irreducibility  of the double conic $\mathscr C_i$ for $i<n-1$, and 
non-reality of each lines of the splitting double conic $\mathscr C_{n-1}$
(see Section \ref{ss:dc}).
Therefore $c_i=c_j$ for any $i,j\in
\{1,2,\cdots,n-1\}$.

Next for the double curves $\mathscr C_n$ and $\mathscr C_{n+1}$,
since $(F|_{H_k\cap Y}) = 2\mathscr C_k
= (Q^2|_{H_k\cap Y})$ for $k\in\{n,n+1\}$
on the cone,
there exists a constant $c_k$ such that 
$F - c_kQ^2\in I_{H_k\cap Y}=(z_{k-1}) + I_Y$.
So taking a difference with $F-c_1Q^2\in I_{ P_1}$,
we obtain that $(c_1-c_k)Q^2 \in (z_{k-1}) + I_Y + 
I_{ P_1}$.
But since $ P_1\subset Y$, we have
$I_{ P_1}\supset I_Y$,
and therefore 
$(c_1-c_k)Q^2 \in (z_{k-1}) + I_{ P_1}$.
Hence if $c_1\neq c_k$ we have $Q^2\in (z_{k-1}) + I_{ P_1}$, which means $Q^2|_{ P_1}\in (z_{k-1}|_{ P_1})$. 
This implies that the divisor $(Q^2)|_{ P_1}$ contains a line $(z_{k-1})$ on 
the plane $ P_1$ as an irreducible component,
which  again contradicts the irreducibility of the double conic $\mathscr C_1$.
Therefore we have $c_1=c_k$ for $k\in \{n,n+1\}$.
By rescaling we can suppose $c_i=1$ 
for any $1\le i\le n+1$.
Thus  we have
\begin{align}
F - Q^2 \in 
\bigg(
\bigcap_{1\le i\le n-1} I_{ P_i}
\bigg)
\cap \big((z_{n-1}) + I_Y\big)
\cap  \big((z_{n}) + I_Y\big).
\label{ideal1}
\end{align}
Let $\Pi$ be the linear subspace of 
$\CP^{n-2}$ spanned by the $(n-2)$ points
$\lmd_1,\cdots,\lmd_{n-2}$ on $\Lmd$.
Since $\Lmd$ is a rational normal curve,   $\Pi$ is $(n-3)$-dimensional.
Let $f\in\CC[z_0,\cdots,z_{n-2}]$ be
a defining linear polynomial of 
the hyperplane $\Pi$.
Then we have 
\begin{align}
 \pi^{-1} (\Pi)  \cap Y = \bigcup_{1\le i\le n-2}  P_{i},
\end{align}
and therefore
$$\bigcap _{1\le i\le n-2} I_{ P_i}
=I_{ P_1\cup\cdots\cup\, P_{n-2}} = I_{\pi^{-1}(\Pi) \, \cap \,Y}
= I_{\pi^{-1}(\Pi)} + I_Y.$$
Hence  \eqref{ideal1} can be rewritten as
\begin{align}
F - Q^2 \in 
\big(
(f) + I_Y
\big)
\cap \big((z_{n-1}) + I_Y\big)
\cap  \big((z_{n}) + I_Y\big)\cap I_{P_{n-1}}.
\label{ideal3}
\end{align}
Further, by an elementary argument,
it is easy to see that 
$$
\big((z_{n-1}) + I_Y\big)
\cap  \big((z_{n}) + I_Y\big)
=
(z_{n-1}z_n) + I_Y
$$
and
$$
\big(
(f) + I_Y
\big)
 \cap
\big((z_{n-1}z_n) + I_Y \big) =
(fz_{n-1}z_n) + I_Y . 
$$
Therefore 
%
%
we can  write
\begin{align}\label{F-Q^2}
F- Q^2 &= z_{n-1} z_n f g + y,\quad y\in I_Y.
\end{align}
Then recalling that this is also in the ideal $I_{P_{n-1}}$,
by restricting both hand sides to the plane 
$P_{n-1}$,
we obtain $g\in I_{P_{n-1}}$.
From our normalization of the homogeneous coordinates 
given just before Theorem \ref{thm:main},
we have $I_{P_{n-1}} = 
(z_0, z_1, \cdots, z_{n-3}).$
Therefore $g\in (z_0,z_1,\cdots,z_{n-3})$.
As $g$ is linear, this means
$g\in \CC[z_0,z_1,\cdots,z_{n-3}]$.
Then if we regard the divisor $(fg)$ as a sum of two hyperplanes
in $\CP^{n-2}$,
since $S_i^++S_i^-$, $1\le i\le n-1$,  are all reducible members of the pencil $|F|$, 
we have, as sets,
\begin{align}\label{intersect3}
(fg)\cap {\Lmd} = \{\lmd_1,\cdots,\lmd_{n-1}\}.
\end{align}
Furthermore, from the definition of the hyperplane $\Pi$,
we have $(f)\cap \Lmd = \{ \lmd_1,\cdots,\lmd_{n-2}\}$,
where all points are included by multiplicity one.
Now if some $\lmd_i$, $1\le i\le n-2$, are contained 
in the intersection $(g)\cap \Lmd$, 
then from  \eqref{F-Q^2},
in a neighborhood of the double conic $\mathscr C_i$,
 the defining equation on $Y$ of the branch divisor $B$ is of the form
$$
q^ 2 + (\lmd-\lmd_i)^k, \quad k\ge 2,
$$
where $\lmd$ is a local coordinate on a neighborhood of $\lmd_i$ in the curve $\Lmd$,
and $q$ is a non-homogeneous representative of $Q$.
This means that the divisor $B$ is singular along
the double conic $\mathscr C_i$, $1\le i\le n-2$.
But if this is actually the case, 
the double cover of $Y$ with branch divisor being $B$ would have singularity along the (inverse image of) 
the real irreducible curve $\mathscr C_i$.
(Note that here we have used the assumption $i\neq n-1$.)
Therefore the morphism $\Phi'_{n-2}:Z'_{n-2}\to Y$ obtained in Section \ref{s:deg1} must contract a real divisor.
But from our explicit elimination, there is no real divisor over the plane $P_i$, $1\le i\le n-2$.
This is a contradiction, and hence we obtain that $(g)\cap \Lambda = \{ \lmd_{n-1}\}$.
This means that $g= cz_0$ for some $c\neq 0$.
Hence, from \eqref{F-Q^2}, we finally obtain
$$
F- Q^2 = z_0 z_{n-1}z_n f+y,
\quad y\in I_Y.
$$
Therefore by disposing $y$,
we obtain the claim of the theorem.
\proofend

\subsection{Dimension of the moduli space}
\label{ss:moduli}
Finally in this subsection we first compute dimension of the moduli space of the present twistor spaces, and then
explain relationship between other twistor spaces.

Let $Z$ be any one of the relevant twistor spaces on $n\CP^2$
and $S$ a real irreducible member of the pencil $|F|$
as before.
Then by a similar argument to \cite[Proposition 5.1]{HonDS4_1},
for the tangent sheaf of $Z$ we have
\begin{align}\label{hfgeg7}
H^i(\Theta_Z) = 0 {\text{ for }} i\neq 1,
\quad h^1(\Theta_Z)=7n-15.
\end{align}
Also, it is easy to show 
\begin{gather}\label{73hew}
H^i(\Theta_S) = 0 {\text{ for }} i\neq 1,
\quad h^1(\Theta_S)=4n-6,\\
h^0(K_S^{-1}) = 1, 
\quad h^1(K_S^{-1})= 2n-8.
\label{839l}\end{gather}
Let $\Theta_{Z,S}$ denote the subsheaf 
of  $\Theta_Z$ consisting 
of germs of vector fields which are tangent to $S$,
and write $\Theta_Z(-S) := \Theta_Z\otimes\mathscr O_Z(-S)$.
Then since $Z$ is Moishezon, we have $H^2(\Theta_Z(-S))
= 0$ by \cite[Lemma 1.9]{C91-2}.
Hence from the cohomology exact sequence
$0 \to \Theta_Z(-S) \to \Theta_{Z,S} \to \Theta_S\to 0$,
by using \eqref{73hew}, we obtain $H^2(\Theta_{Z,S}) = 0$.
Therefore by several standard exact sequences of sheaves including this one, and noting $N_{S/Z}\simeq K_S^{-1}$,
we obtain the following commutative diagram of cohomology groups
on $Z$ and $S$:
$$
\begin{CD}
@. 0 @. 0 @. 0 @.\\
@. @VVV @VVV @VVV @.\\
0 @>>> H^0(\Theta_Z|_S)  @>>> H^0(K_S^{-1})
@>>> H^1(\Theta_S) @>>> H^1(\Theta_Z|_S)\\
@. @VVV @VV{\gamma}V @| @.\\
0 @>>> H^1(\Theta_Z(-S)) @>>> H^1(\Theta_{Z,S})
@>{\alpha}>> H^1(\Theta_S) @>>> 0\\
@. @VVV @VV{\beta}V @VVV @.\\
0 @>>> H^1(\Theta_Z) @= H^1(\Theta_Z) @>>> 0\\
@. @. @VVV\\
@. @. H^1(K_S^{-1})\\
@. @. @VVV\\
@. @. 0
\end{CD}
$$
From the middle column of this diagram we obtain 
$h^1(\Theta_{Z,S}) = 5n-6$ by
\eqref{hfgeg7} and \eqref{839l},
which means $h^1(\Theta_Z(-S)) = n$ from 
the middle row and \eqref{73hew}.

In order to compute the dimension of the moduli space,
we  recall that our twistor spaces can be 
characterized by the property that they have 
one of the rational surface $S$ constructed in 
Section \ref{ss:S} 
as a member of the system $|F|$.
From the construction,
the surface $S$ is determined from 6 points 
on the anticanonical cycle $C_1+C_2+\ol C_1 + \ol C_2$
on $\CP^1\times\CP^1$ (which give the intermediate surface
$S_0$), and therefore  by taking automorphisms of $\CP^1\times\CP^1$
which preserve the cycle into account, they determine 
a 4-dimensional subspace of $H^1(\Theta_S)$.
We denote this subspace  by $V$.
We have $\dim \alpha^{-1}(V) = \dim V + h^1(\Theta_Z(-S))= n+4$.
The tangent space of the moduli space of 
the present twistor spaces can be considered as
the space $\beta(\alpha^{-1}(V))
\subset H^1(\Theta_Z)$.
The image of the map $\gamma$ in the diagram
corresponds to deformations of the pair $(Z,S)$ 
that can be obtained by moving $S$ in $Z$,
and of course they do not give a non-trivial deformation
of $Z$.
Further, from the characterization of $Z$ by
the presence of $S$,
even if we move $S$ in $Z$, the deformed $S$
must still be the one discussed in Section \ref{ss:S}.
This means that the image
of $\gamma$ is contained in $\aaa^{-1}(V)$.
Thus the tangent space  of the moduli space of the present
twistor spaces can be identified with
the quotient space
\begin{align}\label{8fds82}
\aaa^{-1}(V)/\gamma H^0(K_S^{-1}),
\end{align}
and this is $(n+3)$-dimensional by \eqref{839l}.
Thus the dimension of the moduli space is strictly larger
than that of the twistor spaces in \cite{HonDSn1},
which was $n$-dimensional.


Recall that the twistor spaces studied in \cite{HonDSn1} 
also have a structure of a branched double covering over the same scroll $Y$
under the same multiple system $|(n-2)F|$, and the branch divisor 
is a cut of the scroll by a quartic hypersurface.
(In \cite{HonDSn1} the twistor spaces are  presented rather
as a double cover over the resolved space $\tilde Y$,
but it is not difficult to rewrite it as a double cover over the scroll $Y$.)
Hence the structure of the two kinds of the twistor spaces 
is very similar.
Looking defining equations of the quartic hypersurfaces,
or inspecting structure of the surface $S$ in the pencil $|F|$,
it is easy to see that the twistor spaces in \cite{HonDSn1}
are obtained as a limit (under a deformation) of the present twistor spaces.

Next we explain a relationship between the present twistor spaces
(and also those in \cite{HonDSn1}),
and LeBrun twistor spaces \cite{LB91}, from a viewpoint of moduli.
For this, we recall that  
LeBrun's self-dual conformal classes on $n\CP^2$ are determined from distinct $n$ points on 
the hyperbolic space $\mathscr H^3$.
So for each $n\ge 3$, let $H ^{[n]}$ be
the space of configurations of distinct $n$ points on $\mathscr H^3$.
This is a dense open subset of the symmetric product of 
$n$ copies of $\mathscr H^3$.
For each $k$ with $2\le k\le n$, we denote by 
$H ^{[n]}_k$ for the subset of $H ^{[n]}$ consisting of configurations for which the maximal number
of points lying on a common geodesic is exactly $k$.
These provide $H ^{[n]}$ with a natural stratification 
\begin{align}\label{stra}
H ^{[n]}
=H ^{[n]}_2
\supset
H ^{[n]}_3
\supset
H ^{[n]}_4
\supset
\cdots
\supset
H ^{[n]}_{n-1}
\supset
H ^{[n]}_{n}.
\end{align}
The space $H ^{[n]}$ is $3n$-dimensional,
and we clearly have 
$$\dim H ^{[n]}_{k+1} = \dim H ^{[n]}_k -2,
\quad 2\le k<n.$$
In particular we have $\dim H ^{[n]}_k= 3n - 2(k-2)$.
The isometric action of the group ${\rm{PSL}}(2,\CC)$ on
the hyperbolic space
$\mathscr H^3$ naturally induces an action on the space
$H ^{[n]}$ by the same group,
and it clearly preserves the stratification \eqref{stra}.
So we can define the quotients of the strata by
\begin{align}
\mathcal{LB}^{[n]} _k :=  H ^{[n]}_k / {\rm{PSL}}(2,\CC).
\end{align}
The largest space $\mathcal{LB}^{[n]} :=\mathcal{LB}^{[n]} _2$ is exactly
the moduli space of LeBrun's self-dual conformal
classes on $n\CP^2$, and from \eqref{stra} it is equipped with a natural stratification
\begin{align}\label{stra2}
\mathcal{LB}^{[n]} =
\mathcal{LB}^{[n]} _2
\supset
\mathcal{LB}^{[n]} _3
\supset
\mathcal{LB}^{[n]} _4
\supset\cdots\supset
\mathcal{LB}^{[n]} _{n-1}
\supset
\mathcal{LB}^{[n]} _n.
\end{align}
We note that the smallest strata $\mathcal{LB}^{[n]} _n$ is
precisely the  moduli space of toric LeBrun metrics
on $n\CP^2$.
For any $k\neq n$, the ${\rm{PSL}}(2,\CC)$-action on
$H ^{[n]}_k$ is effective and 
therefore we have
\begin{align}\label{dimLB}
\dim\mathcal{LB}^{[n]} _{k} 
= (3n-2k+4) - 6= 3n-2k-2.
\end{align}
On the other hand, for the case $k=n$,
a U(1)-subgroup of ${\rm{PSL}}(2,\CC)$ which fixes the geodesic
acts trivially on 
the smallest stratum $H ^{[n]}_n$, and hence 
we have
\begin{align}\label{LBT}
\dim\mathcal{LB}^{[n]} _n 
= (n+4) - 5= n-1.
\end{align}

This stratification on the moduli space of LeBrun metrics 
is closely related to the moduli space of the present twistor spaces
and also that of the twistor spaces in \cite{HonDSn1} 
in the following way.
As mentioned in the final portion of \cite{HonDSn1},
the twistor spaces in \cite{HonDSn1} 
can be obtained as a small ($\CC^*$-equivariant)
deformation of the twistor spaces of toric LeBrun metrics,
where a $\CC^*$-subgroup of the torus $\CC^*\times\CC^*$ 
is chosen in an (explicit) appropriate way.
Since the moduli space of toric LeBrun metrics on $n\CP^2$ is 
$(n-1)$-dimensional as in \eqref{LBT} and the moduli space
of the twistor spaces on $n\CP^2$ studied in \cite{HonDSn1} is 
$n$-dimensional as computed in the paper, we can conclude that 
the former moduli space is contained in 
the closure of the latter moduli space
as a  hypersurface.
In other words, {\em the moduli space of the twistor spaces
in \cite{HonDSn1} can be partially compactified by attaching the moduli space
of toric LeBrun metrics}, and
{\em the last moduli space is a hypersurface in the partial compactification}.

In order to explain a similar relationship between the present twistor 
spaces and LeBrun twistor spaces on $n\CP^2$, we look at the stratum
\begin{align}\label{n-2}
\mathcal{LB}_{n-2}^{[n]}
\end{align}
in the moduli space of LeBrun metrics on $n\CP^2$.
This is $(n+2)$-dimensional by \eqref{dimLB}.
Now from the construction in Section \ref{ss:S},
it is immediate to see that 
 a member $S\in |F|$ in the present
twistor spaces can be obtained as a small deformation of
a member in $|F|$ on a LeBrun twistor space $Z\in\mathcal{LB}^{[n]}_{n-2}$.
From this, by using a standard argument in deformation theory,
we can show that 
the present twistor spaces are obtained as a small deformation
of a twistor space belonging to $\mathcal{LB}^{[n]}_{n-2}$.
Therefore, recalling that the moduli space of the present
twistor space is $(n+3)$-dimensional as seen above,
we again conclude that 
{\em the moduli space of the twistor spaces studied in this paper
can be partially compactified by attaching the stratum
\eqref{n-2} in the moduli space
of LeBrun metrics on $n\CP^2$}, and
{\em the last moduli space is a hypersurface in the partial compactification}.

One might wonder if a similar relationship between
the moduli spaces carries over
to other strata in the stratification \eqref{stra2}
and twistor spaces of double solid type.
In this respect, 
it seems quite certain that the real situation is  summarized as
in the following diagram:
%
\begin{equation}\label{conj}
 \CD
 \mathcal{LB}\,_{n}^{[n]}@<<< 
 \mathcal{LB}\,_{n-1}^{[n]}@<<<
 \mathcal{LB}\,_{n-2}^{[n]}@<<<
 \mathcal{LB}\,_{n-3}^{[n]}\\
 @AAA @AAA @AAA @AAA\\
 \mathcal{DS}_{\rm{IV}}@<<<
  \mathcal{DS}_{\rm{III}}@<<<
   \mathcal{DS}_{\rm{II}}@<<<
    \mathcal{DS}_{\rm{I}},
 \endCD
 \end{equation}
 where $\mathcal{DS}_{\rm{IV}}$ and 
 $\mathcal{DS}_{\rm{II}}$ are respectively 
 the moduli spaces of the twistor spaces in \cite{HonDSn1}
 and the ones in the present paper, 
 and $\mathcal{DS}_{\rm{III}}$ and 
 $\mathcal{DS}_{\rm{I}}$ are moduli spaces of some unknown
 twistor spaces of double solid type.
Also, in the diagram \eqref{conj}, the notation $A\to B$ means 
that the twistor spaces belonging to the moduli space $B$ are
obtained as a limit (specialization) of 
twistor spaces belonging to the moduli space $A$.
Furthermore, the upper four moduli spaces should always be
a hypersurface in the closure of the lower moduli spaces.
We note that the results in \cite{HonDS4_1, HonDS4_2}
rigorously shows
 this is actually the case when $n=4$.
Actually, the above notations I, II, III and IV for the 
double solid twistor spaces are taken from 
the paper \cite{HonDS4_2}.
We leave the investigation of the `unknown' twistor spaces
belonging to  $\mathcal{DS}_{\rm{I}}$ and 
 $\mathcal{DS}_{\rm{III}}$  in a future paper;
 here we just mention that 
 unlike the ones belonging to $\mathcal{DS}_{\rm{II}}$
 and $\mathcal{DS}_{\rm{IV}}$,
the multiple system $|(n-2)F|$ 
of twistor spaces  belonging to $\mathcal{DS}_I$ is composed with 
a pencil $|F|$, and therefore for analysis of the spaces
we need to consider a linear system of higher degree.
On the other hand, concerning the remaining strata $\mathcal{LB}\,
^{[n]}_k$, $k<n-3$, although we can consider similar
small deformations of LeBrun twistor spaces in these strata,
they are not Moishezon anymore.

The above discussion concerns  relations between
the twistor spaces of double solid type and LeBrun twistor spaces.
Now recall that in \cite{HonLBl} we have constructed 
a family of Moishezon twistor spaces over $n\CP^2$ 
which  share many properties with LeBrun twistor spaces.
If we call these  as LeBrun-like twistor spaces,
 analogous relations to \eqref{conj} hold between
the twistor spaces of double solid type and 
LeBrun-like twistor spaces.
 Namely, we can define a stratification on the moduli space
of LeBrun-like twistor spaces in a similar way
to \eqref{stra} in terms of the structure of a member of 
the system $|F|$, and for the  smallest  four strata among them, 
the relations  \eqref{conj} carry over, 
including the dimension of the moduli spaces.
(We recall that the full moduli space of LeBrun-like twistor spaces
on $n\CP^2$ is $(3n-6)$-dimensional as computed in \cite[Section 5.2]{HonLBl},
which is exactly the same as the original LeBrun metrics.)

Hence it might be possible to say that,
{\em while the moduli spaces of 
LeBrun twistor spaces and LeBrun-like twistor spaces
are not adjacent,
they are interpolated by the moduli space
of the twistor spaces of double solid type.}

%
%
%

\end{document}